# SPECTRAL MOMENTS OF RANDOM MATRICES WITH A RANK-ONE PATTERN OF VARIANCES

VICTOR M. PRECIADO,* *University of Pennsylvania*

M. AMIN RAHIMIAN,* *University of Pennsylvania*

## Abstract

Let $\mathbf{a}_{ij}$, $1 \leq i \leq j \leq n$, be independent random variables and $\mathbf{a}_{ji} = \mathbf{a}_{ij}$, for all $i,j$. Suppose that every $\mathbf{a}_{ij}$ is bounded, has zero mean, and its variance is given by $\sigma_i \sigma_j$, for a given sequence of positive real numbers $\Psi = \{\sigma_i, i \in \mathbb{N}\}$. Hence, the matrix of variances $V_n = (\text{Var}(\mathbf{a}_{ij}))_{i,j=1}^n$ has rank one for all $n$. We show that the empirical spectral distribution of the symmetric random matrix $\mathbf{A}_n(\Psi) = (\mathbf{a}_{ij}/\sqrt{n})_{i,j=1}^n$ converges weakly (and with probability one) to a deterministic limiting spectral distribution which we fully characterize by providing closed-form expressions for its limiting spectral moments in terms of the sequence $\Psi$. Furthermore, we propose a hierarchy of semidefinite programs to compute upper and lower bound on the expected spectral norm of $\mathbf{A}_n$, for both finite $n$ and the limit $n \to \infty$.

Secondary 60B20

## 1. Introduction

We begin by introducing some elementary notation: the set of real and natural numbers are denoted by $\mathbb{R}$ and $\mathbb{N}$, respectively; $\mathbb{N}_0 = \{0\} \cup \mathbb{N}$, $n \in \mathbb{N}$ is a parameter, and $[n]$ denotes $\{1, 2, \ldots, n\}$. The $n \times n$ identity matrix is denoted by $I_n$. Throughout the paper, we will use boldface to indicate random variables and capital letters to denote matrices. Also, an almost sure event is one that occurs with probability one. Consider an $n \times n$ real-valued, symmetric random matrix $\mathbf{A}_n = (\mathbf{a}_{ij}/\sqrt{n})_{i,j=1}^n$ with diagonal and upper triangular entries being independent random variables. Consider further a given sequence of positive real numbers $\Psi = \{\sigma_i \colon i \in \mathbb{N}\}$. We impose the following conditions on our random matrix model $\mathbf{A}_n$ and its characterizing sequence $\Psi$:

> **Assumption 1 (Zero mean).** All entries are zero mean: $\mathbb{E}\{\mathbf{a}_{ij}\} = 0$, for all $i, j$.
>
> **Assumption 2 (Uniformly bounded).** All entries are almost surely bounded in absolute value by a common constant $K > 0$: $\mathbb{P}\{|\mathbf{a}_{ij}| < K\} = 1$, for all $i, j$.
>
> **Assumption 3 (Rank-one variance pattern).** The variances of all entries are specified as follows: $\text{Var}\{\mathbf{a}_{ij}\} = \sigma_i \sigma_j$, for all $i, j$.
>
> **Assumption 4 (Logarithmic growth).** Associated with $\Psi$, there are two monotone sequences $\{\hat{\sigma}_n \colon n \in \mathbb{N}\}$ and $\{\check{\sigma}_n \colon n \in \mathbb{N}\}$ given by $\hat{\sigma}_n = \max_{i \in [n]} \sigma_i$ and $\check{\sigma}_n = \min_{i \in [n]} \sigma_i$. We assume that $\hat{\sigma}_n / \check{\sigma}_n = O(\log n)$.
>
> **Assumption 5 (Mean $k$-th power).** We assume that the $k$-th power means, defined as:
>
> $$\Lambda_k = \lim_{n \to \infty} (1/n) \sum_{i=1}^{n} \sigma_i^k, \qquad (1)$$
>
> exist for all $k \in \mathbb{N}$. Note that Assumptions 2 and 3 together imply that the sequence $\Psi$ is uniformly bounded, so that $\Lambda_k$ is finite (bounded by $K^k$) whenever it exists.

---
* Postal address: Department of Electrical and Systems Engineering, University of Pennsylvania (e-mails: `preciado,mohar@seas.upenn.edu`).

### 1.1. Background and related work

In contrast with most results in the literature [10, 34, 48, 41, 19, 20, 39], the entries of this random matrix ensemble have non-identical variances. With a few exceptions (including the band matrix model described in Appendix C), random matrix ensembles with non-identical variances have only been recently considered [29, 1, 14, 21, 25]. In [21] generalized Wigner matrices with non-identical variances are considered, whose variance matrix is assumed to be doubly stochastic (with unit row- and column-sums). For this ensemble, Wigner's semi-circle still applies and the authors prove universality of the eigenvalue spacing statistics in the bulk. More recent results extend these arguments to random Winger-type matrices with non-identical variances, where the doubly stochasticity condition on the variance matrix is relaxed; and therefore, the semicircle law no longer applies [4]. In this case, the limiting spectral density is a deformation of the semi-circle law that can be characterized through a system of non-linear equations that asymptotically relate the variances to the diagonal entries of the resolvent matrix. The solutions to these quadratic vector equations over the complex upper-half plane are subsequently studied in [2, 3]. As we will show in this paper, a rank-one pattern of variances allows us to apply Wigner's trace method [49, 50] to study the spectral properties of $\mathbf{A}_n$ using combinatorial arguments. Subsequently, we are able to derive closed-form expressions for the limiting spectral moments in terms of the power means $\{\Lambda_k\}_{k \geq 1}$.

Our results are related to both the band matrix model by Anderson and Zeitouni [8] (in Appendix C, we elaborate on this relationship), as well as the free multiplicative convolution (as described in Appendix B). In particular, the limiting spectral distribution of our matrix ensemble with the rank-one variance profile can be also expressed as a free multiplicative convolution between the limiting spectral distribution of a standard Wigner matrix, whose limiting spectral density is semicircular, and the limiting distribution of the entries of the sequence $\Psi = \{\sigma_i, i \in \mathbb{N}\}$. This leads to an alternative derivation of the limiting spectral moments in terms of the limiting $k$-th power means $\{\Lambda_k, k \in \mathbb{N}\}$ using free probability techniques. It is worth highlighting that the free multiplicative convolution techniques are only valid in the limit as the matrix size goes to infinity (for random matrices that are asymptotically free); in contrast, our techniques can be used to bound moments for finite $n$ and give explicit expression for the asymptotic spectral moments as $n \to \infty$. Furthermore, the free probability approach provides only an implicit characterization of the limiting spectral moments. As we discuss in Appendix B, finding expressions for the spectral moments using free probability requires the inversion of a moment generating function that, in general, requires tedious algebraic manipulations (see Appendix B for more details). Finally, whenever the variance profile $\{\sigma_i\}_{i=1}^n$ is estimated from empirical data, the methodology proposed in this paper allows to compute explicit closed-form expression for the spectral moments as a function of the power averages $\Lambda_k$ of the sequence $\{\sigma_i\}_{i=1}^n$, which can be directly computed from the empirical data.

**1.2. Main result**

Before we present the main result in this paper, we need to introduce some nomenclature. Let $\lambda_1(A) \leq \lambda_2(A) \leq \ldots \leq \lambda_n(A)$ be the $n$ real eigenvalues of a symmetric matrix $A$ (ordered from smallest to largest). The spectral radius of $A$ is defined as $\rho(A) = \max_i\{|\lambda_i(A)|\}$. The *empirical spectral measure* of a random symmetric matrix $\mathbf{A}_n$ is defined as $\mathcal{L}_n\{\cdot\} = \frac{1}{n}\sum_{i=1}^{n} \delta_{\lambda_i(\mathbf{A}_n)}\{\cdot\}$, where $\delta_x\{\cdot\}$ is the Dirac delta measure centered at $x$. The *empirical spectral distribution* (ESD) of $\mathbf{A}_n$ is defined as $\mathbf{F}_n(x) = \mathcal{L}_n\{(-\infty, x]\} = \frac{1}{n}\text{card}(\{i \in [n]\colon \lambda_i(\mathbf{A}_n) \leq x\})$, where $\text{card}(\mathcal{X})$ denotes the cardinality of a set $\mathcal{X}$. The $k$-th *empirical spectral moment* of $\mathbf{A}_n$ is defined as $\mathbf{m}_k^{(n)} = \int_{-\infty}^{+\infty} x^k d\mathbf{F}_n(x) = \frac{1}{n}\sum_{i=1}^{n} \lambda_i^k(\mathbf{A}_n)$. The empirical spectral moments are real-valued random variables and their expectations, $\bar{m}_k^{(n)} = \mathbb{E}(\mathbf{m}_k^{(n)})$, are called the *expected spectral moments* of $\mathbf{A}_n$. In this paper, we investigate the limiting behavior of $\mathbf{F}_n(\cdot)$ as $n \to \infty$. Under Assumptions 1 to 5, the random distribution $\mathbf{F}_n(\cdot)$ converges weakly and almost surely (see Appendix A for more details about these modes of convergence) to a deterministic distribution $F(\cdot)$, called the *limiting spectral distribution* (LSD). The *limiting spectral moments* of $\mathbf{A}_n$ are defined as the moments of the LSD, i.e. the $k$-th limiting spectral moment is defined as $m_k = \int_{-\infty}^{+\infty} x^k \, dF(x)$.

Notice that, for the special case $\sigma_i = \sigma$ for all $i \in \mathbb{N}$, we recover the classical Wigner's random matrix ensemble [49, 50]. The LSD in this case is the well-known semi-circle distribution with support $[-2\sigma, 2\sigma]$. The original proof of this classical result was based on computing the limiting spectral moments using the so-called trace method [49, 50]. Accordingly, the $2k$-th spectral moment of the Wigner's random matrix is given by $\frac{\sigma^{2k}}{k+1}\binom{2k}{k}$, which uniquely charactize the semi-circular distribution supported on $[-2\sigma, 2\sigma]$. The first proof of this result relied on the assumption that all entries are identically distributed. This result was later extended by Füredi and Komlós [23] to the case of independent (possibly non-identical) entries with identical variances. The main goal of this paper is to characterize the LSD of random matrices with entries having non-identical variances where the pattern of variances is characterized by a given sequence $\Psi$ satisfying Assumptions 1 to 5. In particular, we provide closed-form expressions for the limiting spectral moments of $\mathbf{A}_n$. We derive these expressions based on a combinatorial argument that allows us to apply the trace method to random matrices with non-identical variances satisfying Assumptions 3 to 5.

It is worth remarking that given an arbitrary sequence $\Psi$, it may not be possible to find an analytic expression for the LSD of the corresponding $\mathbf{A}_n$. Such an analytic expression can only be found for a handful of cases. Consequently, our work does not aim towards finding an analytical expression for the LSD of $\mathbf{A}_n$, but rather to characterize the LSD in terms of the spectral moments. In particular, we provide closed-form expressions for all the limiting spectral moments of the LSD in terms of the power means of the sequence $\Psi$, defined as (1). Our main result is stated as follows.

**Theorem 1.** (Main result.) *The empirical spectral distribution $\mathbf{F}_n(\cdot)$ of the random matrix ensemble $\mathbf{A}_n(\Psi)$ satisfying Assumptions 1 to 5 converges weakly, almost surely, to $F(\cdot)$ as $n \to \infty$, where $F(\cdot)$ is uniquely characterized by the following spectral moments:*

$$m_{2s} = \sum_{(r_1,\ldots,r_s) \in \mathcal{R}_s} \frac{2}{s+1} \binom{s+1}{r_1,\ldots,r_s} \Lambda_1^{r_1} \Lambda_2^{r_2} \ldots \Lambda_s^{r_s}, \text{ and } m_{2s+1} = 0, \qquad (2)$$

*for all $s \in \mathbb{N}_0$, where $\mathcal{R}_s := \{(r_1, r_2, \ldots, r_s) \in \mathbb{N}_0^s \colon \sum_{j=1}^s r_j = s+1, \sum_{j=1}^s j\, r_j = 2s\}$. It is further true that with probability one, $\mathbf{F}_n(\cdot)$ converges weakly to $F(\cdot)$ as $n \to \infty$.*

This result can be used to efficiently compute a truncated sequence of spectral moments of $\mathbf{A}_n$. In practical applications, apart from the spectral moments, one usually cares about other spectral properties of $\mathbf{A}_n$. For example, the spectral radius finds applications in many practical scenarios, including combinatorics, mathematical physics, and theoretical computer science. In Section 3, we propose a hierarchy of semidefinite programs (SDP) to compute upper and lower bounds on the spectral radius of $\mathbf{A}_n$ given a finite sequence of spectral moments.

## 2. Proof of the main result

We prove our main result using an extension of Wigner's trace method. In its standard form, Wigner's method is used to prove the semicircle law for random matrices with i.i.d. entries presenting identical variances (and satisfying certain technical conditions). Wigner's proof is based on a combinatorial argument that relates the spectral moments with the moments of the semicircular distribution. In the analysis of random matrices with non-identical variances, this combinatorial argument is not directly applicable. In this paper, we refine Wigner's combinatorial argument to derive closed-form expressions for the spectral moments of random matrices whose entries have non-identical variances following a rank-one pattern (Assumption 3).

### 2.1. Moments and convergence

Before we introduce the combinatorial elements of our proof, we describe the convergence properties of our results. Choose an element $\mathbf{i}$ from the set $[n] = \{1, \ldots, n\}$ uniformly at random. Let $\lambda(\mathbf{A}_n) = \lambda_{\mathbf{i}}(\mathbf{A}_n)$ be the random variable that corresponds to the $\mathbf{i}$-th real eigenvalue of $\mathbf{A}_n$. Denote the law of the random variable $\lambda(\mathbf{A}_n)$ by $\overline{\mathcal{L}}_n\{\cdot\}$. Consider the *expected spectral distribution* that is associated with this law: $\overline{F}_n(x) = \overline{\mathcal{L}}_n\{(-\infty, x]\} = \mathbb{P}\{\lambda(\mathbf{A}_n) \leq x\}$, for all $x \in \mathbb{R}$. The relation between $\overline{F}_n(\cdot)$ and the empirical spectral distribution $\mathbf{F}_n(\cdot)$ is better understood upon noting that for all $x \in \mathbb{R}$, $\overline{F}_n(x) = \mathbb{P}\{\lambda(\mathbf{A}_n) \leq x\} = \mathbb{E}\{\mathcal{L}_n\{(-\infty, x]\}\} = \mathbb{E}\{\mathbf{F}_n(x)\}$, where the second equality follows from Fubini-Tonelli theorem. Similarly, for any $k, n \in \mathbb{N}$, the $k$-th expected spectral moment satisfies

$$\bar{m}_k^{(n)} = \int_{-\infty}^{+\infty} x^k d\overline{F}_n(x) = \mathbb{E}\left\{\int_{-\infty}^{+\infty} x^k d\mathbf{F}_n(x)\right\}.$$

Applying the definition of $\mathbf{F}_n(x)$ to the right-hand side, we get

$$\bar{m}_k^{(n)} = \mathbb{E}\left\{\frac{1}{n}\sum_{i=1}^{n}\lambda_i^k(\mathbf{A}_n)\right\} = \mathbb{E}\left\{\frac{1}{n}\text{trace}(\mathbf{A}_n^k)\right\} = \mathbb{E}\{\mathbf{m}_k^{(n)}\}, \qquad (3)$$

where $\mathbf{m}_k^{(n)}$ is the $k$-th spectral moment of the random matrix $\mathbf{A}_n$. In writing (3), we have implicitly presumed almost sure finiteness of $\mathbf{m}_k^{(n)}$, as well as integrability of $x^k$ against $d\overline{F}_n(x)$; these facts are verified in Appendix D for the random matrix ensemble satisfying Assumptions 1-5.

Our proof of Theorem 1 is based on the method of moments, which characterizes the limiting spectral distribution from the values of the expected spectral moments: $\lim_{n\to\infty} \bar{m}_k^{(n)} = m_k$, for all $k \in \mathbb{N}$. As stated below (and proved in Appendix D), for our random matrix ensemble $\mathbf{A}_n(\Psi)$, pointwise convergence of the expected spectral moments implies almost sure and weak convergence of the ESD $\mathbf{F}_n(\cdot)$:

**Theorem 2.** (The method of moments.) *If $\lim_{n\to\infty} \bar{m}_k^{(n)} = m_k$ for all $k \in \mathbb{N}$, then the ESD $\mathbf{F}_n(\cdot)$ of the random matrix $\mathbf{A}_n$ satisfying Assumptions 1-5 converges weakly,* almost surely, to the LSD $F(\cdot)$ as $n \to \infty$, where $F(\cdot)$ is the unique distribution function satisfying $\int_{-\infty}^{+\infty} x^k \, dF(x) = m_k$ for all $k \in \mathbb{N}$. It is further true that, with probability one, $\mathbf{F}_n(\cdot)$ converges weakly to $F(\cdot)$ as $n \to \infty$.

Our proof of the Theorem 1 proceeds by first showing that the expected spectral distributions $\overline{F}_n(\cdot)$ converge weakly to the distribution $F(\cdot)$ as $n \to \infty$, where $F(\cdot)$ is the unique distribution with the moments sequence $\{m_k \colon k \in \mathbb{N}\}$, which are given by (2). For our proof, we follow the method of moments and the technicalities are spelled out in Appendix D. The results of Appendix D, as summarized by Theorem 2, allow us to conclude the almost sure and weak convergence of the ESD $\mathbf{F}_n(\cdot)$ to the deterministic distribution $F(\cdot)$, from the pointwise convergence of the expected spectral moments $\{\bar{m}_k^{(n)} \colon k \in \mathbb{N}\}$ to the sequence of moments $\{m_k \colon k \in \mathbb{N}\}$, as $n \to \infty$.

The convergence proof for the moments is executed in two steps. First, we identify those terms that asymptotically dominate the behavior of each spectral moment (in Subsection 2.3). Second, we derive the asymptotically exact expressions for each of the dominant terms (in Subsection 2.4). Our main techniques are based on a refinement of the combinatorial argument used to prove Wigner's semicircle law. The required combinatorial preliminaries are established next.

## 2.2. Combinatorial preliminaries

The Catalan numbers admit many combinatorial interpretations and play a central role in the combinatorial proof of Wigner's semicircle law [7]. An integer sequence $\{b_t\}_{t=0}^{2s}$ such that $b_0 = b_{2s} = 0, b_t \geq 0$, and $|b_{t+1} - b_t| = 1$ for all $t = 0, 1, 2, \cdots, 2s-1$ is called a *Dyck path* of length $2s$. The $s$-th Catalan number, defined as $C_s = \frac{1}{s+1}\binom{2s}{s}$, counts the total number of Dyck paths of length $2s$, cf. [7]. Catalan numbers can also be used to count *non-crossing partitions* of an ordered set [7], as well as many other combinatorial structures [35]. Specially relevant to our work is the relationship between the Catalan numbers and *rooted ordered trees*. A rooted ordered tree $T$ is a tree in which one vertex is designated as the root and the children of each vertex are ordered (see [35], page 221); i.e. there is a total order $\precsim$ on the vertex set of $T$, respecting the partial order $\precsim$ defined as follows: for all $\{j, k\} \subset \mathcal{V}(T)$, $j \precsim k$ iff $j$ belongs to the unique path on $T$ that connects $k$ to the root. It is possible to construct a bijection between Dyck paths of length $2s$ and ordered trees with $s$ edges (see [7, Lemma 2.1.6

]); hence, the number of ordered trees with $s$ edges is equal to $C_s$. Although these combinatorial entities play an important role in the proof of the semicircle law, they are not enough to study random matrices whose entries have non-identical variances.

Given a connected graph $G$ with $s$ edges and $n$ vertices labeled $\{1,\ldots,n\}$, we denote by $d_i(G)$ the degree of vertex $i$ (i.e., the number of edges incident to $i$ in $G$). The degree distribution of $G$ is defined as the sequence of integers $\overline{r}(G) = (r_1,\ldots,r_s) \in \mathbb{N}_0^s$, where $r_d = r_d(G)$ is the number of vertices with degree $d$ in $G$. We drop the graph argument $(G)$ when there is no danger of confusion: using $d_i$ for the degree of vertex $i$ and $r_d$ for number of vertices with degree $d$. Notice that the maximum degree is at most $s$; hence, $r_d = 0$ for all $d > s$. Given a connected graph $G$ with $s$ edges, $G$ is a tree if and only if it has $s+1$ vertices [35]. Hence, we can easily characterize the set of degree distributions corresponding to trees, as follows:

**Proposition 1.** (Degree distribution of trees.) *Consider a connected graph $G$ with $s$ edges and degree distribution $(r_1,\ldots,r_s) \in \mathbb{N}_0^s$. Then, $G$ is a tree if and only if $\sum_{d=1}^{s} r_d(G) = s+1$ and $\sum_{d=1}^{s} d\, r_d = 2s$.*

*Proof.* The proof is based on the following simple observations. First, for any graph $G$ with $s$ edges, its degree distribution $(r_1,\ldots,r_s)$ satisfies $\sum_{d=1}^{s} d\, r_d = 2s$. Furthermore, the summation $\sum_{d=1}^{s} r_d = s+1$ implies that $G$ has $s+1$ vertices. The statement in the proposition is therefore true, since a connected graph with $s$ edges and $s+1$ vertices is always a tree. □

**Remark 1.** Based on Proposition 1, notice that the set $\mathcal{R}_s$ defined in Theorem 1 denotes the set of integer sequences that are valid degree distributions for trees with $s$ edges.

For $s \in \mathbb{N}$, let $\mathcal{T}_s$ denote the set of all ordered trees on $s$ vertices. Given a degree distribution $\overline{r} \in \mathbb{N}_0^{s-1}$, the set $\mathcal{T}_s(\overline{r})$ is defined as the subset of ordered trees in $\mathcal{T}_s$ with degree distribution $\overline{r}$. A classical result in enumerative combinatorics [35] states that the total number of ordered trees on $s+1$ vertices is given by $\operatorname{card}(\mathcal{T}_{s+1}) = C_s$, where $C_s$ is the $s$-th Catalan number, defined by $C_s = \frac{1}{s+1}\binom{2s}{s}$. Furthermore, the following lemma follows almost directly from Theorem 5.3.10 in [35]:

**Lemma 1.** (Counting rooted ordered trees.) *For $s \in \mathbb{N}$, let $\overline{r} = (r_1,...,r_s) \in \mathbb{N}_0^{s+1}$,*

with $\sum_{j=1}^{s} r_j = s+1$ and $\sum_{j=1}^{s} j\, r_j = 2s$. Then, the number of rooted ordered trees with $s+1$ nodes presenting a degree distribution equal to $\bar{r}$ is given by

$$\operatorname{card}\left(\mathcal{T}_{s+1}(\bar{r})\right) = \frac{2}{s+1}\binom{s+1}{r_1, r_2, \ldots, r_s}.$$

**Remark 2.** An interesting Catalan identity can be directly obtained from this result. Indeed, using the partition $\mathcal{T}_{s+1} = \bigcup_{\bar{r}\in\mathcal{R}_s} \mathcal{T}_{s+1}(\bar{r})$, we can express the total number of ordered trees on $s+1$ vertices as, $C_s = \operatorname{card}(\mathcal{T}_{s+1}) = \sum_{\bar{r}\in\mathcal{R}_s} \operatorname{card}(\mathcal{T}_{s+1}(\bar{r}))$. Hence, from Lemma 1, we obtain

$$C_s = \frac{2}{s+1} \sum_{\bar{r}\in\mathcal{R}_s} \binom{s+1}{r_1, r_2, \ldots, r_s}. \tag{4}$$

Indeed, if we let $\sigma_i = 1$ for all $i$, then it follows that $\Lambda_k = 1, \forall k$, and by (4) and after replacing in (2), we recover the moment sequence for the classical Wigner semi-circle law: $m_{2s} = C_{2s}$ and $m_{2s-1} = 0$ for all $s \in \mathbb{N}$ (see [49, 50]).

Lemma 1 plays an important role in determining closed-form expressions of the limiting spectral moments in (2). As we will show below, there is a correspondence between rooted ordered trees and a particular class of closed walks on the complete graph. We will prove that this class of closed walks asymptotically dominate the expression of the moments. To show this, let us Consider the complete graph $\mathcal{K}_n$ with vertices labeled by $[n]$. A closed walk $w$ of length $k$ on $\mathcal{K}_n$ is an ordered finite sequence of integers $w = (i_1, i_2, \ldots, i_{k-1}, i_k, i_1)$ such that $i_j \in [n]$ for all $j$. Define the sets of vertices and edges visited by $w$ as $\mathcal{V}(w) = \{i_j \colon j \in [k]\}$ and $\mathcal{E}(w) = \{\{i_j, i_{j+1}\}\colon j \in [k-1]\}\cup\{\{i_k, i_1\}\}$. For any $e \in \mathcal{E}(w)$, we define $N(e, w)$ as the number of times that walk $w$ transverses the edge $e$ in any direction. We denote by $\mathcal{W}_k$ the set of all closed walks of length $k$ on $\mathcal{K}_n$. It is useful to partition the set $\mathcal{W}_k$ into subsets $\mathcal{W}_{k,p}$ defined as $\mathcal{W}_{k,p} = \{w \in \mathcal{W}_k \colon \operatorname{card}(\mathcal{V}(w)) = p\}$, i.e. the set of closed walks of length $k$ visiting $p$ vertices. Furthermore, it is convenient to define the following subset of $\mathcal{W}_{k,p}$:

$$\widehat{\mathcal{W}}_{k,p} = \{w \in \mathcal{W}_{k,p} \colon N(e, w) \geq 2 \text{ for all } e \in \mathcal{E}(w)\}, \tag{5}$$

i.e., the set of walks in $\mathcal{W}_{k,p}$ for which each edge is traversed at least twice.

A key step in Wigner's trace method is to count the number of walks in $\widehat{\mathcal{W}}_{2s, s+1}$,

which is given by Lemma 2.3.15 in [37]:

$$\mathrm{card}(\widehat{\mathcal{W}}_{2s,s+1}) = n(n-1)\ldots(n-s)C_s. \qquad (6)$$

The main idea behind the proof of this identity is to establish a bijection between the set of walks in $\widehat{\mathcal{W}}_{2s,s+1}$ and the set of rooted ordered trees (or, alternatively, the set of non-crossing cycles of length $k$; see Section 2.3.5 in [37] for more details). We begin by observing that the subgraph induced by the set of edges visited by any walk in $\widehat{\mathcal{W}}_{2s,s+1}$ is always a tree with $s+1$ vertices. This is true because of the following facts: $(i)$ on $s+1$ vertices, we need at least $s$ distinct edges for the induced graph to be connected, and $(ii)$ for any walk in $\widehat{\mathcal{W}}_{k,p}$, each edge is visited at least twice, which implies that for any walk in $\widehat{\mathcal{W}}_{2s,s+1}$ each edge must be visited *exactly* twice; hence, the total number of edges in the induced graph is exactly $s$, which implies that this graph is a tree. Furthermore, any walk $w \in \widehat{\mathcal{W}}_{2s,s+1}$ induces a total order on the vertices according to the order of the first appearance of each vertex in the walk $w$. Therefore, a rooted order tree uniquely encodes a walk $w \in \widehat{\mathcal{W}}_{2s,s+1}$, as follows: $(i)$ the $s$ edges of the tree indicate the $s$ edges visited by the walk $w$ of length $2s$, each edge being visited twice (each time in a different direction); $(ii)$ the root of the tree encodes the initial vertex of the walk; and $(iii)$ the total order of vertices in the tree encodes the order in which the vertices are visited along the walk $w$ (similar to a depth-first traversal of a tree, as described in [17]).

### 2.3. Wigner's trace method

In this subsection, we briefly review Wigner's trace method. The ideas in this subsection are well known and we simply adapt them to the case of random matrices with a rank-one pattern of variances. Wigner's trace method is based on the following simple observation about the trace of the powers of $\mathbf{A}_n = (\mathbf{a}_{ij}/\sqrt{n})_{i,j=1}^n$:

$$\bar{m}_k^{(n)} = \mathbb{E}\left\{\frac{1}{n}\mathrm{trace}(\mathbf{A}_n^k)\right\} = \frac{1}{n^{k/2+1}} \sum_{1 \leq i_1, i_2, \ldots, i_k \leq n} \mathbb{E}\left\{\mathbf{a}_{i_1 i_2}\mathbf{a}_{i_2 i_3}\cdots\mathbf{a}_{i_{k-1} i_k}\mathbf{a}_{i_k i_1}\right\}. \quad (7)$$

To each closed walk of length $k$ in the complete graph $\mathcal{K}_n$, denoted by $w = (i_1, i_2, \ldots, i_{k-1}, i_k, i_1)$, we associate a weight $\omega(w) = \mathbb{E}\left\{\mathbf{a}_{i_1 i_2}\mathbf{a}_{i_2 i_3}\cdots\mathbf{a}_{i_{k-1} i_k}\mathbf{a}_{i_k i_1}\right\}$; hence, (7) can be

rewritten as

$$\bar{m}_k^{(n)} = \frac{1}{n^{k/2+1}} \sum_{w \in \mathcal{W}_k} \omega(w) = \frac{1}{n^{k/2+1}} \sum_{p=1}^{k} \sum_{w \in \mathcal{W}_{k,p}} \omega(w) = \sum_{p=1}^{k} \mu_{k,p}, \quad (8)$$

where we have used the partition $\mathcal{W}_k = \cup_{p \in [k]} \mathcal{W}_{k,p}$ in the second equality, and $\mu_{k,p}$ in the last equality is defined by

$$\mu_{k,p} = \frac{1}{n^{k/2+1}} \sum_{w \in \mathcal{W}_{k,p}} \omega(w). \quad (9)$$

The latter can be simplified upon noting that only a subset of walks in $\mathcal{W}_{k,p}$ has a non-zero contribution to the summation in (9). In particular, the set of walks presenting a non-zero weight is $\widehat{\mathcal{W}}_{k,p}$, defined in (5). This is true since for any $w \in \mathcal{W}_{k,p} \setminus \widehat{\mathcal{W}}_{k,p}$, there exists an edge $\{i,j\} \in \mathcal{E}(w)$ such that $N(\{i,j\}, w) = 1$. Hence, by independence,

$$\omega(w) = \mathbb{E}\{\mathbf{a}_{ij}\} \mathbb{E}\left\{ \prod_{\{k,l\} \in \mathcal{E}(w) \setminus \{\{i,j\}\}} \mathbf{a}_{kl}^{N(\{k,l\},w)} \right\} = 0,$$

since $\mathbb{E}\{\mathbf{a}_{ij}\} = 0$. Therefore, $\mu_{k,p}$ in (9) simplifies into

$$\mu_{k,p} = \frac{1}{n^{k/2+1}} \sum_{w \in \widehat{\mathcal{W}}_{k,p}} \omega(w). \quad (10)$$

An important observation in the proof of Wigner's semicircle law [49, 50] is that for $k$ even, the summation in (8) is asymptotically dominated by those closed walks belonging to $\widehat{\mathcal{W}}_{k,k/2+1}$; i.e. the term $\mu_{k,k/2+1}$. In what follows, we find conditions under which the same kind of dominance holds for the random matrix ensemble $\mathbf{A}_n$ satisfying Assumptions 1-5. We proceed with two lemmas, the first of which follows by a trite counting argument, while the second is at the heart of several existing results [23, 15, 48].

**Lemma 2.** *It holds true for all $p > \lfloor k/2 \rfloor + 1$ that $\mu_{k,p} = 0$; wherefore, (8) simplifies to $\bar{m}_k^{(n)} = \sum_{p=1}^{\lfloor k/2 \rfloor + 1} \mu_{k,p}$, for all $k, n \in \mathbb{N}$.*

*Proof.* It follows from the pigeonhole principle that every walk $w \in \mathcal{W}_{k,p}$ with $p > \lfloor k/2 \rfloor + 1$ has at least one edge $e \in \mathcal{E}(w)$ such that $N(e, w) = 1$, whence $\widehat{\mathcal{W}}_{k,p} = \varnothing$ for $p > \lfloor k/2 \rfloor + 1$, and from (10) we get that $\mu_{k,p} = 0$, as claimed. □

**Lemma 3.** *It holds true for all $p \in [n]$ that*

$$|\mu_{k,p}| \leq \frac{1}{n^{k/2+1-p}} \binom{k}{2p-2} p^{2(k-2p+2)} 4^{p-1} K^{k-2p+2} \hat{\sigma}_n^{2p-2}.$$

*Proof.* First, note that for all $w \in \widehat{\mathcal{W}}_{k,p}$ and all $\{i,j\} \in \mathcal{E}(w)$, $N(\{i,j\},w) \geq 2$, while $|\mathbf{a}_{ij}| < K$, almost surely, so that $\left|\mathbb{E}\left\{\mathbf{a}_{ij}^{N(\{i,j\},w)}\right\}\right| \leq \mathbb{E}\left\{\left|\mathbf{a}_{ij}^{N(\{i,j\},w)}\right|\right\} \leq \hat{\sigma}_n^2 K^{N(\{i,j\},w)-2}$. The latter together with the independence of the random entries $\mathbf{a}_{ij}$ imply that for all $w \in \widehat{\mathcal{W}}_{k,p}$, $|\omega(w)| \leq K^{k-2p+2}\hat{\sigma}_n^{2p-2}$, where we have used the fact that with $p$ distinct vertices in walk $w$, there is at least $p-1$ distinct edges. For each one of these edges, we can use the bound $\hat{\sigma}_n^2$ leading to the $\hat{\sigma}_n^{2(p-1)}$ term, which we multiply by $K^{k-2p+2}$ to account for those edges whose multiplicities are greater than 2. Next, we make use of the following bound (which is developed in Section 3.2 of [23] and is subsequently used in Lemma 2 of [15] and Equation (5) of [48]): for all $p \in [n]$

$$\operatorname{card}\left(\widehat{\mathcal{W}}_{k,p}\right) \leq n(n-1)\ldots(n-p+1)\binom{k}{2p-2}p^{2(k-2p+2)}2^{2p-2}$$

$$\leq n^p \binom{k}{2p-2}p^{2(k-2p+2)}2^{2p-2}. \tag{11}$$

The claim now follows from (10) and upon combining (11) with the bound, $|\omega(w)| \leq K^{k-2p+2}\hat{\sigma}_n^{2p-2}$, derived above for all $w \in \widehat{\mathcal{W}}_{k,p}$. □

We now proceed to our main dominance result; in particular, we prove that for $k = 2s$ and under certain conditions the term $\mu_{2s,s+1}$ dominates the other terms of the summation in (8).

**Theorem 3.** (Dominant walks for even moments.) *If* $\frac{s^6}{\hat{\sigma}_n^2}\left(\frac{\hat{\sigma}_n}{\check{\sigma}_n}\right)^{2s} = o(n)$, *then* $\bar{m}_{2s}^{(n)} = (1+o(1))\mu_{2s,s+1}$.

*Proof.* Note per Lemma 2 that $\bar{m}_{2s}^{(n)} = \sum_{p=1}^{s+1} \mu_{2s,p}$. We show the desired dominance by first lower bounding the term $\mu_{2s,s+1}$ and then upper bounding the terms $\mu_{2s,p}, p < s+1$ as follows. We begin by noting from (6) that

$$\operatorname{card}\left(\widehat{\mathcal{W}}_{2s,s+1}\right) = n(n-1)\ldots(n-s)\frac{1}{s+1}\binom{2s}{s} \geq \frac{(n-s)^{s+1}}{s+1}\binom{2s}{s}. \tag{12}$$

We next lower-bound the contribution of each walk in $\widehat{\mathcal{W}}_{2s,s+1}$ as $\omega(w) \geq \check{\sigma}_n^{2s}$, which holds true because $N(e,w) = 2, \forall e \in \mathcal{E}(w)$, and together with (12) implies that

$$\mu_{2s,s+1} \geq \left(\frac{n-s}{n}\right)^{s+1}\frac{1}{s+1}\binom{2s}{s}\check{\sigma}_n^{2s}. \tag{13}$$

Now, for the case $p < s+1$, we apply Lemma 3 with $k = 2s$ to get

$$\mu_{2s,p} \leq \frac{1}{n^{s+1-p}}4^{p-1}\binom{2s}{2p-2}p^{4(s-p+1)}\hat{\sigma}_n^{2p-2}K^{2s-2p+2}. \tag{14}$$

We next form the ratio between the two inequalities (13) and (14) to get

$$\begin{aligned}\frac{\mu_{2s,s+1}}{\mu_{2s,p}} &\geq \frac{1}{2}\left(\frac{4}{K^2p^4s^2}\right)^{s+1-p}\frac{\check{\sigma}_n^{2s}(n-s)^{s+1}}{\hat{\sigma}_n^{2p-2}n^p}\\ &\geq \frac{2\hat{\sigma}_n^2}{K^2s^6}\left(\frac{\check{\sigma}_n}{\hat{\sigma}_n}\right)^{2s}\left(\frac{n-s}{n}\right)^s(n-s),\end{aligned} \quad (15)$$

where in the first inequality we have used that

$$\binom{2s}{2p-2} \leq \left(\frac{2s}{2s-2p+2}\right)^{2s-2p+2} \leq 2s^{2s-2p+2},$$

and in the second inequality we take into account that the greatest lower bound is achieved for $p = s$, cf. [23]. The proof now follows upon noting that if $\frac{s^6}{\hat{\sigma}_n^2}\left(\frac{\hat{\sigma}_n}{\check{\sigma}_n}\right)^{2s} = o(n)$, then $\mu_{2s,p} = o(1)\mu_{2s,s+1}$; thence, $\bar{m}_{2s}^{(n)} = (1+o(1))\mu_{2s,s+1}$, as claimed. □

**Theorem 4.** (Vanishing odd moments.) *Under Assumptions 1 to 5, we have that $\bar{m}_{2s+1}^{(n)} = o(1)$.*

*Proof.* First note, by Lemma 3, that each of the terms $\mu_{2s+1,p}, p \leq s+1$ can be upper bounded as follows

$$\begin{aligned}|\mu_{2s+1,p}| &\leq \frac{n^{(p-1)}}{n^{s+(1/2)}}\binom{2s+1}{2p-2}p^{2(2s-2p+3)}2^{2p-2}K^{2s-2p+3}\hat{\sigma}_n^{2p-2}\\ &\leq \frac{1}{n^{s+(1/2)}}\binom{2s+1}{s}(s+1)^{2(2s+3)}2^{2s}K^{2s+1}\left(\frac{n\hat{\sigma}_n^2}{K^2}\right)^{p-1}.\end{aligned} \quad (16)$$

Writing $\bar{m}_{2s+1}^{(n)} = \sum_{p=1}^{s+1}\mu_{2s+1,p}$ by Lemma 2 and replacing from (16), the odd spectral moments can now be upper bounded as

$$\begin{aligned}|\bar{m}_{2s+1}^{(n)}| &\leq \frac{1}{n^{s+(1/2)}}\binom{2s+1}{s}(s+1)^{2(2s+3)}2^{2s}K^{2s+1}\sum_{p=1}^{s+1}\left(\frac{n\hat{\sigma}_n^2}{K^2}\right)^{p-1}\\ &= \frac{1}{n^{s+(1/2)}}\binom{2s+1}{s}(s+1)^{2(2s+3)}2^{2s}K^{2s+1}\frac{\left(n\left(\frac{\hat{\sigma}_n}{K}\right)^2\right)^{s+1}-1}{n\left(\frac{\hat{\sigma}_n}{K}\right)^2-1},\end{aligned}$$

and it follows that, if $\hat{\sigma}_n^{2s} = o(\sqrt{n})$, then $\bar{m}_{2s+1}^{(n)} = o(1)$. Given the bounded entries and the specified variances, we get that under Assumptions 1 to 5, $\hat{\sigma}_n^{2s}$ is bounded uniformly in $n$ and the claim follows. □

### 2.4. Counting dominant walks

In the following theorem, we use the tools from Section 2.2 to derive explicit expressions for the dominant term $\mu_{2s,s+1}$, which will allow us to write closed-form expressions for the limiting spectral moments.

**Theorem 5.** (Asymptotically exact expressions for even-order moments.) *For all $s \in \mathbb{N}$, if $\hat{\sigma}_n^{2s} = o(n)$, then*

$$\lim_{n \to \infty} \mu_{2s,s+1} = \sum_{\bar{r} \in \mathcal{R}_s} \frac{2}{s+1} \binom{s+1}{r_1, \ldots, r_s} \prod_{j=1}^{s} \Lambda_j^{r_j} = m_{2s},$$

*where $m_{2s}$ is defined in (2).*

*Proof.* Starting from the definition (10), we have that

$$\mu_{2s,s+1} = \frac{1}{n^{s+1}} \sum_{w \in \widehat{\mathcal{W}}_{2s,s+1}} \omega(w)$$

$$= \frac{1}{n^{s+1}} \sum_{T \in \mathcal{T}_{s+1}} \prod_{\{i,j\} \in \mathcal{E}(T)} \sigma_i \sigma_j,$$

where in the last equality we take into account the fact that $\mathbb{E}\{\mathbf{a}_{ij}^2\} = \sigma_i \sigma_j$ and rewrite the summation over the set of walks $\widehat{\mathcal{W}}_{2s,s+1}$ in terms of the set of underlying rooted ordered trees; here, again, we are relying critically on the bijection between the set of walks in $\widehat{\mathcal{W}}_{2s,s+1}$ and the set of rooted ordered trees on $s$ vertices (please, refer to (6) and the explanations therein for more details). For any tree $T$, it is always true that $\prod_{\{i,j\} \in \mathcal{E}(T)} \sigma_i \sigma_j = \prod_{i \in \mathcal{V}(T)} \sigma_i^{d_i(T)}$; where we use $d_i(T)$ for the degree of node $i$ in the tree $T$; hence, the above can be written as:

$$\mu_{2s,s+1} = \frac{1}{n^{s+1}} \sum_{T \in \mathcal{T}_{s+1}} \left\{ \sum_{\substack{1 \leq i_1, \ldots, i_{s+1} \leq n \\ \text{card}(\{i_1, \ldots, i_{s+1}\}) = s+1}} \left( \prod_{k=1}^{s+1} \sigma_{i_k}^{d_k(T)} \right) \right\}, \quad (17)$$

where the second summation is over the choice of $s+1$ distinct vertices out of the set $[n]$ for each tree $T \in \mathcal{T}_{s+1}$. To simplify the preceding expression, consider any rooted ordered tree $T$ on the set of nodes $[s+1]$ with degrees given by the sequence $(d_1(T), \ldots, d_{s+1}(T))$, whose degree distribution is $\bar{r} := (r_1(T), \ldots, r_s(T))$. Define the $d$-th power average of the sequence $\{\sigma_k \colon k \in [n]\}$ as $S_{n,d} = \sum_{k=1}^{n} \sigma_k^d$. Using the definition of degree distribution $\bar{r}$ and its relation to the degrees $d_j(T)$ of different nodes $j \in [s+1]$, we can write:

$$\prod_{j=1}^{s} S_{n,j}^{r_j(T)} = \prod_{j=1}^{s+1} S_{n,d_j(T)} = \prod_{j=1}^{s+1} \left( \sum_{k=1}^{n} \sigma_k^{d_j(T)} \right)$$

$$= \sum_{i_1=1}^{n} \sum_{i_2=1}^{n} \cdots \sum_{i_{s+1}=1}^{n} \sigma_{i_1}^{d_1(T)} \sigma_{i_2}^{d_2(T)} \cdots \sigma_{i_{s+1}}^{d_{s+1}(T)}$$

$$= \sum_{j=1}^{s+1} \left\{ \sum_{\substack{1 \leq i_1, \ldots, i_{s+1} \leq n \\ \text{card}(\{i_1, \ldots, i_{s+1}\}) = j}} \left( \prod_{k=1}^{s+1} \sigma_{i_k}^{d_k(T)} \right) \right\}. \qquad (18)$$

Next, we have that for all $j \in [s]$

$$\lim_{n \to \infty} \frac{1}{n^{s+1}} \sum_{j=1}^{s} \left\{ \sum_{\substack{1 \leq i_1, \ldots, i_{s+1} \leq n \\ \text{card}(\{i_1, \ldots, i_{s+1}\}) = j}} \left( \prod_{k=1}^{s+1} \sigma_{i_k}^{d_k(T)} \right) \right\} = 0, \qquad (19)$$

which is true since the term in the curly brackets above can be bounded as

$$0 \leq \frac{1}{n^{s+1}} \sum_{\substack{1 \leq i_1, \ldots, i_{s+1} \leq n \\ \text{card}(\{i_1, \ldots, i_{s+1}\}) = j}} \left( \prod_{k=1}^{s+1} \sigma_{i_k}^{d_k(T)} \right) \leq \frac{1}{n^{s+1}} \binom{n}{j} \hat{\sigma}_n^{2s} \leq \frac{n^{j-s-1}}{j!} \hat{\sigma}_n^{2s}. \qquad (20)$$

Next, note that Assumptions 2 and 3 imply that $\hat{\sigma}_n^{2s}$ is bounded uniformly in $n$; hence, the right-hand side of the inequalities in (20) is $o(1)$. By taking the limits of both sides in (17), and combining (18) and (19), we get

$$\lim_{n \to \infty} \mu_{2s,s+1} = \sum_{T \in \mathcal{T}_{s+1}} \lim_{n \to \infty} \frac{1}{n^{s+1}} \sum_{\substack{1 \leq i_1, \ldots, i_{s+1} \leq n \\ \text{card}(\{i_1, \ldots, i_{s+1}\}) = s+1}} \prod_{k=1}^{s+1} \sigma_{i_k}^{d_k(T)}$$

$$= \sum_{T \in \mathcal{T}_{s+1}} \lim_{n \to \infty} \frac{1}{n^{s+1}} \prod_{j=1}^{s} S_{n,j}^{r_j(T)}$$

$$= \sum_{T \in \mathcal{T}_{s+1}} \prod_{j=1}^{s} \Lambda_j^{r_j(T)},$$

where in the last equality we use that $\Lambda_j = \lim_{n \to \infty} (1/n) S_{n,j}$. To finish the proof, use the partition $\mathcal{T}_{s+1} = \bigcup_{\bar{r} \in \mathcal{R}_s} \mathcal{T}_{s+1}(\bar{r})$ to get:

$$\lim_{n \to \infty} \mu_{2s,s+1} = \sum_{\bar{r} \in \mathcal{R}_s} \sum_{T \in \mathcal{T}_{s+1}(\bar{r})} \prod_{j=1}^{s} \Lambda_j^{r_j(T)}$$

$$= \sum_{\bar{r} \in \mathcal{R}_s} \operatorname{card}\left(\mathcal{T}_{s+1}(\bar{r})\right) \prod_{j=1}^{s} \Lambda_j^{r_j}$$

$$= \sum_{\bar{r} \in \mathcal{R}_s} \frac{2}{s+1} \binom{s+1}{r_1, \ldots, r_s} \prod_{j=1}^{s} \Lambda_j^{r_j},$$

where Lemma 1 is invoked in the last equality, concluding the proof. $\square$

This brings us to end of the proof for Theorem 1 as Theorems 3, 4 and 5 together imply that $\lim_{n \to \infty} \bar{m}_k^{(n)} = m_k, \forall k \in \mathbb{N}$. Hence, by the method of moments (Theorem 2) we can claim the weak and almost sure convergence of the ESD $\mathbf{F}_n(\cdot)$ as stated in Theorem 1.

### 3. Bounds on the spectral radius using the spectral moments

In Theorem 1, we provide closed-form expressions for the spectral moments of the limiting spectral distribution (LSD) of the random matrix ensemble $\mathbf{A}_n(\Psi)$. In Subsection 3.1, we build on our previous combinatorial analysis to bound the expected spectral moments of finite matrices and compute upper and lower bounds on the expected spectral radius of $\mathbf{A}_n$, for any finite $n$. Next, in Subsection 3.2, we propose a semidefinite program to compute optimal lower bounds on the almost-sure deterministic limit of the spectral radius, $\lim_{n \to \infty} \rho(\mathbf{A}_n)$, given a truncated sequence of spectral moments, $\{m_{2s}, s = 1, \ldots, \hat{s}\}$, which can be efficiently computed from from (2).

### 3.1. Bounding the expected spectral radius for finite random matrices

In this subsection, we present methods to bound the expected spectral radius $\mathbb{E}\{\rho(\mathbf{A}_n)\}$ for any finite $n \in \mathbb{N}$ by using the following well-known inequalities (see Equation (2.66) of [37]),

$$\left(\bar{m}_{2s}^{(n)}\right)^{1/2s} \leq \mathbb{E}\{\rho(\mathbf{A}_n)\} \leq \left(n \bar{m}_{2s}^{(n)}\right)^{1/2s}. \tag{21}$$

To this end, since the results of Theorem 1 are asymptotic in nature and applicable only as $n \to \infty$, we need to refine the combinatorial analysis of Subsection 2.4 to compute upper and lower bounds on the expected spectral moments $\bar{m}_{2s}^{(n)}$ for finite $n$.

**Lemma 4.** (Lower bounds on the expected spectral moments.) *For any finite $n \in \mathbb{N}$, we have that*

$$\bar{m}_{2s}^{(n)} \geq \check{m}_{2s}^{(n)} = \frac{1}{n^{s+1}} \sum_{(r_1,\ldots,r_s) \in \mathcal{R}_s} \frac{2}{s+1} \binom{s+1}{r_1,\ldots,r_s} S_{n,1}^{r_1} \ldots S_{n,s}^{r_s} - \check{\varepsilon}_{2s}^{(n)},$$

*where $\check{\varepsilon}_{2s}^{(n)} = \frac{\hat{\sigma}_n^{2s}}{n^{s+1}} \sum_{j=1}^{s} \binom{n}{j} j^{s+1-j}$.*

*Proof.* From Lemma 2, we get that $\bar{m}_{2s}^{(n)} = \sum_{p=1}^{s+1} \mu_{2s,p} \geq \mu_{2s,s+1}$ and by (17) and (18), we know that $\mu_{2s,s+1}$ is given by

$$\frac{1}{n^{s+1}} \sum_{(r_1,\ldots,r_s) \in \mathcal{R}_s} \frac{2}{s+1} \binom{s+1}{r_1,\ldots,r_s} S_{n,1}^{r_1} \ldots S_{n,s}^{r_s} - \frac{1}{n^{s+1}} \sum_{j=1}^{s} \sum_{\substack{1 \leq i_1,\ldots,i_{s+1} \leq n \\ \text{card}(\{i_1,\ldots,i_s\})=j}} \prod_{k=1}^{s+1} \sigma_{i_k}^{d_k(T)}.$$

Furthermore, we can lower bound the above expression by substituting $\hat{\sigma}_n$ for $\sigma_{i_k}$ to get

$$\bar{m}_{2s}^{(n)} \geq \frac{1}{n^{s+1}} \sum_{(r_1,\ldots,r_s) \in \mathcal{R}_s} \frac{2}{s+1} \binom{s+1}{r_1,\ldots,r_s} S_{n,1}^{r_1} \ldots S_{n,s}^{r_s} - \frac{1}{n^{s+1}} \sum_{j=1}^{s} \sum_{\substack{1 \leq i_1,\ldots,i_{s+1} \leq n \\ \text{card}(\{i_1,\ldots,i_s\})=j}} \hat{\sigma}_n^{2s}$$

$$= \frac{1}{n^{s+1}} \sum_{(r_1,\ldots,r_s) \in \mathcal{R}_s} \frac{2}{s+1} \binom{s+1}{r_1,\ldots,r_s} S_{n,1}^{r_1} \ldots S_{n,s}^{r_s} - \frac{1}{n^{s+1}} \sum_{j=1}^{s} j^{s+1-j} \binom{n}{j} \hat{\sigma}_n^{2s},$$

where we have used the fact that $\sum_{k=1}^{s+1} d_k(T) = 2s$ in the first inequality, and the last expression follows by a trite combinatorial argument, leading to the claimed result. □

The following lower bound on the expected spectral radius is then immediate upon applying Lemma 4 to (21).

**Corollary 1.** (Lower bound on the expected spectral radius.) *For any finite $n$, we have that*

$$\mathbb{E}\{\rho(\mathbf{A}_n)\} \geq n^{-\frac{s+1}{2s}} \left( \sum_{(r_1,\ldots,r_s) \in \mathcal{R}_s} \frac{2}{s+1} \binom{s+1}{r_1,\ldots,r_s} S_{n,1}^{r_1} \ldots S_{n,s}^{r_s} - \check{\varepsilon}_{2s}^{(n)} \right)^{1/2s}.$$

We next compute an upper bound on the expected spectral moments and, consequently, the expected spectral radius.

TABLE 1: Even-order moments of $\mathbf{A}_{1000}$ with $\sigma_i = e^{-4i/n}, \forall i \in [1000]$

| $2s$ | $m_{2s}$ | $\check{m}_{2s}^{(1000)}$ | $\hat{m}_{2s}^{(1000)}$ | Mean (100x) |
|---|---|---|---|---|
| 4  | 1.599e-2 | 1.394e-2 | 2.152e-2 | 1.591e-2 |
| 6  | 5.739e-3 | 4.733e-3 | 5.740e-3 | 5.700e-3 |
| 8  | 2.388e-3 | 1.958e-3 | 2.388e-3 | 2.368e-3 |
| 10 | 1.081e-3 | 0.897e-3 | 1.081e-3 | 1.071e-3 |

**Lemma 5.** (Upper bounds on the expected spectral moments.) *For any finite $n$, we have that $\bar{m}_{2s}^{(n)} \leq \hat{m}_{2s}^{(n)} = (1 + \theta_{2s}^{(n)} s)m_{2s}$, where $m_{2s}$ is defined in (2) and*

$$\theta_{2s}^{(n)} = (K^2 s^6/2n\hat{\sigma}_n^2)(\hat{\sigma}_n/\check{\sigma}_n)^{2s}(n^{s-1})/(n-s)^{s+1}. \tag{22}$$

*Proof.* Recall from the proof of Theorem 5 that $\mu_{2s,s+1}$ approaches its limit $m_{2s}$ from below and, by Lemma 2, we get $\bar{m}_{2s}^{(n)} = \mu_{2s,s+1} + \sum_{p=1}^{s} \mu_{2s,p} \leq m_{2s} + \sum_{p=1}^{s} \mu_{2s,p}$. The claim follows after invoking (15) to upper bound each of the terms $\mu_{2s,p}, \forall p \in [s]$ as $\mu_{2s,p} \leq \theta_{2s}^{(n)} \mu_{2s,s+1}$. □

We can now combine the result of Lemma 5 with the right-hand side inequality in (21) to upper bound the expected spectral radius as follows.

**Corollary 2.** (Upper bound on the expected spectral radius.) *For any finite $n$, we have that $\mathbb{E}\{\rho(\mathbf{A}_n)\} \leq \left(n(1 + \theta_{2s}^{(n)} s)m_{2s}\right)^{1/2s}$, where $\theta_{2s}^{(n)}$ is defined in (22).*

The above bounds are next verified in numerical simulations. We consider a random matrix with $n = 1000$ and a sequence $\sigma_i = e^{-4i/n}$ for $i \in [1000]$. Fig. 1 depicts the histogram of the eigenvalues of one realization of $\mathbf{A}_{1000}$. Notice that, as the entries have non-identical variances, the observed distribution departs significantly from the classical Wigner's semicircle law. In Table 1, we compare the value of the asymptotic spectral moments in Theorem 1 with the lower bounds $\check{m}_{2s}^{(1000)}$ and upper bound $\hat{m}_{2s}^{(1000)}$ on the expected spectral moments for $n = 1000$, as derived in Lemmas 4 and 5, respectively. We also include in this Table the empirical mean values for the spectral moments obtained from averaging 100 realizations of $\mathbf{A}_{1000}$. Notice how the empirical averages (last column in Table 1), as well as the upper bounds $\hat{m}_{2s}^{(1000)}$ are very close to our theoretical predictions for the even-order spectral moments $m_{2s}$.

We next compare the empirical mean of the spectral radius of $\mathbf{A}_{1000}$ for 100 re-alizations with the upper and lower bounds in Corollaries 1 and 2, respectively. The empirical mean is computed as 0.7679, while the lower and upper bounds in Corollaries 1 and 2 for $s = 30$ and $n = 1000$ are 0.6913 and 0.7757, respectively. Tighter lower bounds can be computed in the asymptotic case using the methodology presented in the next subsection.

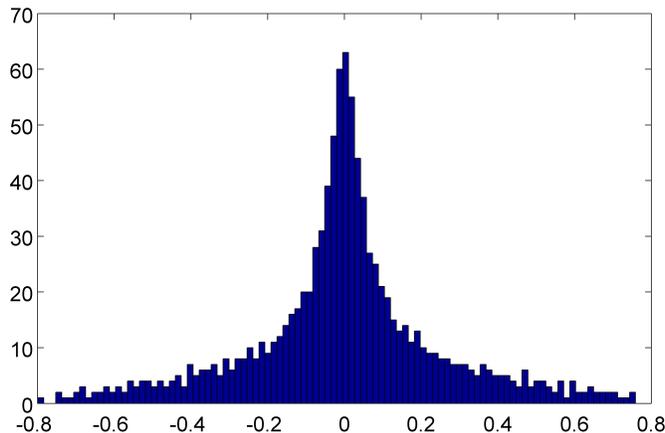

FIGURE 1: The histogram of the eigenvalues of a random realization of $\mathbf{A}_{1000}$ with $\sigma_i = e^{-4i/n}, \forall i \in [1000]$.

### 3.2. Almost sure bounds on the asymptotic spectral radius

We now shift our focus to bounds on the asymptotic spectral radius of $\mathbf{A}_n$ that hold almost surely. To begin, note by Fatou's lemma that $\mathbb{E}\left\{\liminf_{n\to\infty} \rho(\mathbf{A}_n)\right\} \leq \liminf_{n\to\infty} \mathbb{E}\left\{\rho(\mathbf{A}_n)\right\}$. Moreover, from Corollary 4 and the almost sure deterministic limit set forth in Theorem 1, we get that $\lim_{n\to\infty} \rho(\mathbf{A}_n) = \mathbb{E}\{\liminf_{n\to\infty} \rho(\mathbf{A}_n)\}$, almost surely; in addition, as a consequence of Corollary 2, we have $\liminf_{n\to\infty} \mathbb{E}\{\rho(\mathbf{A}_n)\} \leq \liminf_{n\to\infty} \left(n(1+\theta_{2s}^{(n)}s)m_{2s}\right)^{1/2s}$. The preceding facts yield an almost-sure upper bound on the asymptotic spectral radius as

$$\lim_{n\to\infty} \rho(\mathbf{A}_n) \leq \liminf_{n\to\infty} \left(n(1+\theta_{2s}^{(n)}s)m_{2s}\right)^{1/2s},$$

which holds true for any $s \in \mathbb{N}$. Indeed, letting $s \to \infty$ and taking $n = s$ yields

$$\lim_{n\to\infty} \rho(\mathbf{A}_n) \leq \lim_{s\to\infty} m_{2s}^{1/2s},$$

almost surely. In what follows, we use the methodology introduced in [30] to find almost sure lower bounds on $\lim_{n\to\infty} \rho(\mathbf{A}_n)$ that are optimal given the knowledge of a truncated sequence of asymptotic spectral moments from (2), as $\{m_{2s}, s = 1, \ldots, \hat{s}\}$ for $\hat{s}$ a fixed odd integer.

First, note that since $m_{2s+1} = 0$ for all $s \geq 0$, the asymptotic spectral distribution $F(\cdot)$ is symmetric. Associated with $F(\cdot)$, we define the auxiliary distribution $G(\cdot)$ given by $G(x) = F(\sqrt{x})$ for $x \geq 0$ and $G(x) = 0$ when $x < 0$. Denote the $s$-th moment of $G(\cdot)$ by $\nu_s$; hence, we have that $\nu_s = \int_{-\infty}^{+\infty} x^s \, dG(x) = \int_{-\infty}^{+\infty} x^{2s} \, dF(x) = m_{2s}$ for any $s \in \mathbb{N}$. Given a truncated sequence of moments $(\nu_0, \nu_1, \ldots, \nu_{2\bar{s}+1})$ with $\nu_0 = 1$ and $\bar{s} = (\hat{s} - 1)/2$, we can use Proposition 1 in [30] to compute a lower bound on $\sup_{x>0}\{G(x) > 0\}$ by solving the following semidefinite program [45]

$$B_s = \min_{x>0} \quad x$$
$$\text{subject to} \quad \begin{bmatrix} \nu_0 & \nu_1 & \cdots & \nu_{\bar{s}} \\ \nu_1 & \nu_2 & \cdots & \nu_{\bar{s}+1} \\ \vdots & \vdots & \ddots & \vdots \\ \nu_{\bar{s}} & \nu_{\bar{s}+1} & \cdots & \nu_{2\bar{s}} \end{bmatrix} x - \begin{bmatrix} \nu_1 & \nu_2 & \cdots & \nu_{\bar{s}+1} \\ \nu_2 & \nu_3 & \cdots & \nu_{\bar{s}+2} \\ \vdots & \vdots & \ddots & \vdots \\ \nu_{\bar{s}+1} & \nu_{\bar{s}+2} & \cdots & \nu_{2\bar{s}+1} \end{bmatrix} \succeq 0, \quad (23)$$

where $\succeq 0$ indicates the belonging relation to the convex cone of all real symmetric positive semidefinite $(\bar{s} + 1) \times (\bar{s} + 1)$ matrices. Notice that the matrices involved in (23) present a Hankel structure and $B_s$ is the solution of a semidefinite program in one variable $x > 0$, which can be solved efficiently using off-the-shelf software tools. From $B_s$, we get that $\lim_{n\to\infty} \rho(\mathbf{A}_n) \geq \sqrt{B_s}$, almost surely. We verify this bound numerically for $\bar{s} = 14$ to find that the optimal lower bound for the asymptotic spectral radius is 0.7578, which is comparable to the empirical mean 0.7679 computed in Subsection 3.1 for 100 realizations of $\mathbf{A}_{1000}$ with $\sigma_i = e^{-4i/n}, \forall i \in [1000]$. Notice how this latter bound is substantially tighter than the value of 0.6913 obtained from Corollary 1.

## 4. Conclusions

In this paper, we have analyzed a random matrix ensemble characterized by independent, zero-mean, and uniformly bounded entries presenting a rank-one pattern of variances, i.e., $\text{Var}\{\mathbf{a}_{ij}\} = \sigma_i \sigma_j$, for a given sequence $\{\sigma_i : i \in \mathbb{N}\}$. Our main result, stated in Theorem 1, establishes that the spectral distribution of the random matrix ensemble $\mathbf{A}_n(\Psi)$ satisfying Assumptions 1 to 5 converges almost surely and weakly to a deterministic distribution that can be characterized via its spectral moments, for which we provide closed-form expressions. We would like to remark that, even though the spectral moments can be implicitly characterized using the band matrix model or free probability, these alternative approaches require tedious algebraic manipulations to obtain explicit expressions for the moments. In contrast, Theorem 1 directly provides an explicit expression for the asymptotic spectral moments. Based on our analysis, we have also provided upper and lower bounds on the expected spectral moments of *finite-dimensional* random matrices with zero-mean entries and the specified pattern of variances. These bounds can be used to, for example, compute bounds on the expected spectral radius of $\mathbf{A}_n(\Psi)$, for a finite $n$. Finally, we have illustrated how the exact asymptotic expressions for the moments of the LSD can be used to optimally bound the almost sure limit of the spectral radius of $\mathbf{A}_n(\Psi)$ as $n \to \infty$ using semidefinite programming.

### Appendix A. Weak convergence of the empirical spectral distributions

Let $F_n(\cdot)$ and $F(\cdot)$ be distribution functions, then by definition $F_n(\cdot)$ converges weakly to $F(\cdot)$ if $\lim_{n\to\infty} F_n(x) = F(x)$ for each $x \in \mathbb{R}$ at which $F(\cdot)$ is continuous (see Section 14 of [13]). Alternatively, we can use the following as an equivalent characterization.

**Lemma 6.** (Weak convergence, Theorem 4.4.2 of [16].) *Given distribution functions $F_n(\cdot)$ and $F(\cdot)$, $F_n(\cdot)$ converges weakly to $F(\cdot)$ if, and only if, for any bounded continuous real-valued function $f(\cdot)$ we have that*

$$\lim_{n\to\infty} \int_{-\infty}^{+\infty} f(x)\, dF_n(x) = \int_{-\infty}^{+\infty} f(x)\, dF(x).$$

In the case of empirical distributions, such as $\mathbf{F}_n(\cdot)$, we can leverage Lemma 6 to derive the following two variations of weak convergence in probabilistic settings.

**Definition 1.** (*Almost sure weak convergence of the ESD.*) The empirical distribution $\mathbf{F}_n(\cdot)$ is said to converge weakly to $F(\cdot)$, with probability one, if there is a measurable set $\mathcal{S}$ such that $\mathbb{P}\{\mathcal{S}\} = 1$ and for all sample points $\omega \in \mathcal{S}$, the realization $\mathbf{F}_n(\cdot) \mid_\omega$ converges weakly to $F(\cdot)$.

**Definition 2.** (*Weak almost sure convergence of the ESD.*) The empirical distribution $\mathbf{F}_n(\cdot)$ is said to converge weakly, almost surely, to $F(\cdot)$, if for any bounded continuous real-valued function $f(\cdot)$ we have that $\int_{-\infty}^{+\infty} f(x) \, d\mathbf{F}_n(x)$ converges to $\int_{-\infty}^{+\infty} f(x) \, dF(x)$, almost surely.

## Appendix B. Relations to the free multiplicative convolution

The discipline of free probability was founded by Voiculescu from his study of non-commutative random variables in the 1980s; it has subsequently found major applications in the study of high-dimensional random matrices which are non-commutative objects. In particular, free probability identifies a certain sufficient condition of asymptotic freeness under which the asymptotic spectrum of the sums or products can be analyzed using the individual asymptotic spectra of the summands or multiplicands, based on the respective concepts of free additive and multiplicative convolutions. Given a Hermitian random matrix ensemble $A = \{\mathbf{A}_n, n \in \mathbb{N}\}$, let $\phi(A)$ be its limiting expected spectral moment of the first-order, given by $\phi(A) = \lim_{n \to \infty} (1/n) \mathbb{E}\{\text{trace}(A)\}$. Two Hermitian random matrix ensembles $A$ and $B$ are asymptotically free (or freely independent) if for all $l$ and for all polynomials $p_i(\cdot)$ and $q_i(\cdot)$, $1 \leq i \leq l$ satisfying $\phi(p_i(A)) = \phi(q_i(B)) = 0$ for all $i$, we have that $\phi(p_1(A)q_2(B) \ldots p_l(A)q_l(B)) = 0$ [7, Chapter 5]. Accordingly, if $A$ and $B$ are freely independent with the respective limiting spectral distributions $\mu$ and $\nu$, then the free multiplicative convolution of $\mu$ and $\nu$ which is denoted by $\mu \boxtimes \nu$ and it is the limiting spectral distribution of $B^{1/2} A B^{1/2}$ or equivalently $A^{1/2} B A^{1/2}$ [7, Definition 5.3.28]. In this Appendix we will show how the the desired limiting spectral distribution $F$ which is the focus of our main result in Theorem 1 can be expressed as a free multiplicative convolution between the semicircle law and the limiting distribution of the entries of the variance profile sequence $\{\sigma_i, i \in \mathbb{N}\}$. Next we use transform techniques, in particular $S$-transforms, to express the $S$-transform of $F$ (defined below) as the product of the $S$-transforms two distributions involved in the free multiplicative convolution. Finally, by relating the $S$-transforms to the moments of the underlying distributions, we are able to verify the explicit expressions of moments in Theorem 1, indeed satisfy the free multiplicative convolution relation claimed here in the Appendix.

Consider an arbitrary distribution $\mu$ with moments $\{m_k^\mu = \int x^k \mu(dx)\}_{k \geq 0}$ and let

$$\psi_\mu(z) = \sum_{k \geq 1} m_k^\mu z^k = M_\mu(z) - 1.$$

where $M_\mu(z)$ is the complex-valued moment-generating function of $\mu$. Let $\chi_\mu(\cdot)$ be the functional inverse of $\psi_\mu(\cdot)$ satisfying

$$\psi_\mu(\chi_\mu(z)) = z, \qquad (25)$$

The $S$-transform of $\mu$, denote by $S_\mu$, is another complex-valued function which can be computed to the moments and the inverse function $\chi_\mu(\cdot)$ via the following relationships [7, Definition 5.3.28]:

$$S_\mu(z) = \frac{1+z}{z} \chi_\mu(z). \qquad (26)$$

This relation is more useful in the equivalent format $\psi_\mu(z) S_\mu(\psi_\mu(z)) = z(1 + \psi_\mu(z))$ [11, Section 10.2, Equation (17)], which can be used to relate the $S$-transform $S_\mu(z)$ directly to the moments $\{m_k^\mu\}_{k \geq 0}$ rather than through the inverse function $\chi_\mu(z)$. Let us expand the $S$-transform as:

$$S_\mu(z) = \sum_{k \geq 1} s_k^\mu z^{k-1}.$$

The sequence of real numbers $\{s_k^\mu\}_{k \geq 0}$ (i.e., the coefficients of the $S$-transform) can be deduced from the moment sequence $\{m_k^\mu\}_{k \geq 0}$ (and vice versa) using the following identities [11, Section 10.2, Equation (18)]:

$$m_1^\mu s_1^\mu = 1, \qquad (27)$$

$$m_K^\mu = \sum_{k=1}^{K+1} s_k^\mu \sum_{\substack{l_1, \ldots, l_k \geq 1 \\ l_1 + \ldots + l_k = K+1}} m_{l_1}^\mu \cdots m_{l_k}^\mu. \qquad (28)$$

We are now equipped with the necessary tools to study the limiting spectral distribution $F$ as a free multiplicative convolution between the semicircle law and the limiting distribution of the entries of $\Psi = \{\sigma_i, i \in \mathbb{N}\}$. Consider an $n$-dimensional real symmetric standard Wigner matrix $\mathbf{W}_n = \mathbf{W}_n^T$, i.e. with i.i.d. zero-mean, unit-variance entries (in the upper triangular). Given the infinite sequence $\{\sigma_i\}_{i\geq 1}$, also consider the diagonal matrix $\Sigma_n = \text{diag}(\sigma_1, \ldots, \sigma_n)$. Define the matrix $\widehat{\mathbf{A}}_n = \frac{1}{\sqrt{n}}\Sigma_n^{1/2}\mathbf{W}_n\Sigma^{1/2}$, whose $i,j$-th entry has zero mean and variance $\sigma_i\sigma_j/n$, which is the same as the variance profile of the random matrix ensemble $\mathbf{A}_n(\Psi)$ studied in our main result. Indeed, the random matrix ensemble $\widehat{\mathbf{A}}_n$ includes $\mathbf{A}_n(\Psi)$ as a special case, and establishing the limiting spectral distribution of $\widehat{\mathbf{A}}_n$ provides an alternative characterization of the desired eigenvalue distribution $F$, which is the subject of our main result (Theorem 1). Indeed, our objective in this appendix is to characterize the limiting spectral distribution of $\widehat{\mathbf{A}}_n$, or equivalently $\widetilde{\mathbf{A}}_n = \frac{1}{\sqrt{n}}\Sigma_n\mathbf{W}_n$ (both having the same eigenvalue spectrum), as a function of the given sequence $\{\sigma_i\}_{i\geq 1}$.

At the crux of our analysis is the fact that a standard Wigner matrix and a deterministic diagonal matrix (in particular, the matrices $\mathbf{W}_n$ and $\Sigma_n$) are asymptotically free [44, Example 2.38]; Voiculescu [47] proved this fact for the special case of Gaussian entries in 1991, and it was later extended by Dykema (1993) to include all Wigner matrices and deterministic block diagonal matrices as well [18]. In what follows, we build on these results to analyze the limiting spectral distribution of $\widetilde{\mathbf{A}}_n$ using techniques from free probability theory, and in particular the method of free multiplicative convolution. Let us denote by $\theta$ and $\nu_{sc}$ the limiting spectral distributions of $\Sigma_n$ and $\mathbf{W}_n$, respectively. In particular, since $\mathbf{W}_n$ is a Wigner matrix, $\nu_{sc}$ follows Wigner's semicircle law [49, 50] with the following semicircular density:

$$f_{sc}(x) = \frac{1}{2\pi}\sqrt{4-x^2},$$

for $x \in [-2, 2]$ and $f_{sc}(x) = 0$, otherwise. Also, the limiting spectral distribution of the diagonal matrix $\Sigma_n$ is given by $\theta\{\cdot\} = \lim_{n\to\infty} \frac{1}{n}\sum_{i=1}^{n} \delta_{\sigma_i}\{\cdot\}$; here the limit is in the topology of weak convergence and it is assumed existent for the free probability techniques to apply. Using the free multiplicative convolution the limiting spectral distribution of $\mathbf{A}_n$, denoted by $F$, satisfies $F = \theta \boxtimes \nu_{sc}$. Our next step is to characterize $F$ using its moment-generating function and $S$-transform, denoted by $M_F(\cdot) = 1+\psi_F(\cdot)$ and $S_F$ respectively (as defined above). An important property of the $S$-transform is the fact that given two (asymptotically) free matrices, such as $\Sigma_n$ and $\mathbf{W}_n$, the $S$-transform of the (limiting) spectral distribution of the product (i.e., $F = \theta \boxtimes \nu_{sc}$) is given by [7, Lemma 5.3.30]:

$$S_F = S_{sc}S_\theta, \tag{29}$$

where $S_{sc} = 1/\sqrt{z}$ is the $S$-transform of the semi-circular [31, Section 3], which plays a crucial role in our following derivations. The relation (29) was first developed by Voiculescu for random variables with non-zero mean (or non-negative support) [46]. Rao and Speicher [31] extend the applicability of this relation to include the cases where one of the variables has a vanishing mean as in the case of a semicircular density which is symmetric around the origin with zero mean. The extension relies critically on the fact that the inverse function $\chi_\mu$ in (26) can written as a formal power series in $\sqrt{z}$. Using this result, we now derive expressions for the moments of $F$ as a function of the given sequence $\{\sigma_i, i \in \mathbb{N}\}$. We first characterize $S_\theta$ from the sequence $\{\sigma_i, i \in \mathbb{N}\}$, as follows. Notice that the moments of $\theta\{\cdot\}$ are given by

$$m_k^\theta = \lim_{n\to\infty} \frac{1}{n}\sum_{i=1}^n \sigma_i^k = \Lambda_k \tag{30}$$

for all $k \in \mathbb{N}$, which are well-defined per Assumption 5. Hence, we can compute the coefficients $\{s_k^\theta\}_{k\geq 0}$ of $S_\theta$ using (27) and (28), as follows. For example, (27) gives the first coefficient $s_1^\theta$ as:

$$s_1^\theta = \frac{1}{m_1^\theta} = \frac{1}{\Lambda_1}. \tag{31}$$

The second coefficient $s_2^\theta$ can be computed from (28) with $K=1$, as follows:

$$\Lambda_1 = s_1^\theta \Lambda_2 + s_2^\theta \Lambda_1^2 = \frac{\Lambda_2}{\Lambda_1} + s_2^\theta \Lambda_1^2,$$

which implies

$$s_2^\theta = \frac{\Lambda_1^2 - \Lambda_2}{\Lambda_1^3}. \tag{32}$$

Similarly, using (28) with $K=2$, we obtain

$$\Lambda_2 = s_1^\theta \Lambda_3 + 2s_2^\theta \Lambda_1 \Lambda_2 + s_3^\theta \Lambda_1^3$$
$$= \frac{\Lambda_3}{\Lambda_1} + 2\frac{\Lambda_1^2 \Lambda_2 - \Lambda_2^2}{\Lambda_1^2} + s_3^\theta \Lambda_1^3,$$

which implies

$$s_3^\theta = \frac{\Lambda_2}{\Lambda_1^3} - \frac{\Lambda_3}{\Lambda_1^4} - 2\frac{\Lambda_1^2 \Lambda_2 - \Lambda_2^2}{\Lambda_1^5}. \tag{33}$$

Following this methodology, it would be possible, in principle, to compute higher-order coefficients $s_k^\theta$ for each $k \geq 1$ as a function of the sequence of power means $\{\Lambda_j\}_{j=1}^k$, defined in (30). In the sequel, we shall use the sequence of coefficients $\{s_1^\theta, \ldots, s_k^\theta\}$ to compute a truncated sequence of moments of $F$. To this end, we can apply (29) with $S_{sc}(z) = 1/\sqrt{z}$ and $S_\theta(z) = \sum_{k\geq 1} s_k^\theta z^{k-1}$ to obtain:

$$S_F(z) = \frac{1}{\sqrt{z}} \left( \sum_{k\geq 1} s_k^\theta z^{k-1} \right) = \frac{s_1^\theta}{\sqrt{z}} + s_2^\theta \sqrt{z} + s_3^\theta \sqrt{z}^3 + \ldots$$

From (25), we can compute the inverse function of the moment-generating function as:

$$\chi_F(z) = \frac{z}{1+z} S_F(z) = \frac{1}{1+z} \left( s_1^\theta \sqrt{z} + s_2^\theta \sqrt{z}^3 + s_3^\theta \sqrt{z}^5 + \ldots \right).$$

Since $\frac{1}{1+z} = 1 - z + z^2 - z^3 + \ldots$, we have that

$$\chi_F(z) = \left(1 - z + z^2 - z^3 + \ldots\right) \left( s_1^\theta \sqrt{z} + s_2^\theta \sqrt{z}^3 + s_3^\theta \sqrt{z}^5 + \ldots \right)$$
$$= s_1^\theta \sqrt{z} + \left(s_2^\theta - s_1^\theta\right) \sqrt{z}^3 + \left(s_3^\theta - s_2^\theta + s_1^\theta\right) \sqrt{z}^5 - \ldots$$

Therefore, we can write

$$\chi_F(z) = \sum_{k=1}^\infty \beta_k \sqrt{z}^k, \quad \text{where } \beta_k = \begin{cases} \sum_{i=1}^{\frac{k-1}{2}} (-1)^{\frac{k-1}{2}-i} s_i^\theta, & \text{for } k \text{ odd,} \\ 0, & \text{for } k \text{ even.} \end{cases} \quad (34)$$

Next note that the expansion $\psi_F(z) = \sum_{k\geq 1} m_k z^k$ can be combined with the identity $\psi_F(\chi_F(z)) = z$ from (25) to relate the moments sequence $\{m_k\}_{k\geq 1}$ with the coefficients $\{\beta_k\}_{k\geq 1}$. This relationship is given in [31, Proof of Proposition 2.3] as follows:

$$1 = (\beta_1)^2 m_2, \quad (35)$$

$$0 = \sum_{k=2}^r \sum_{\substack{1\leq l_1,\ldots,l_k \leq r \\ l_1+\cdots+l_k=r}} m_k \beta_{l_1} \cdots \beta_{l_k}, \text{ for } r > 2. \quad (36)$$

Coup de grâce is to replace for the coefficients $\{\beta_k\}_{k\geq 1}$ in terms of $\{s_k^\theta\}_{k\geq 1}$ using (34), and then write the latter in terms of $\{\Lambda_k\}_{k\geq 1}$ (as in (31), (32) and (33)) to finally express the moments $\{m_k\}_{k\geq 1}$ in terms of the power means $\{\Lambda_k\}_{k\geq 1}$. Starting with (35), $\beta_1 = s_1^\theta$ in (34), and (31), we get $1 = m_2 \left(s_1^\theta\right)^2 = m_2/\Lambda_1^2$; therefore, $m_2 = \Lambda_1^2$, indeed verifying (2) for $s = 1$. Similarly, for $r = 3$, (36) implies that $0 = 2m_2\beta_1\beta_2 + m_3\left(\beta_1\right)^3 = \frac{m_3}{\Lambda_1^3}$, where in the last equality we have used: $\beta_1 = 1/\Lambda_1$ and $\beta_2 = 0$. Therefore, we have that $m_3 = 0$ as expected from (2). For $r = 4$, (36) with $\beta_1 = 1/\Lambda_1$, $\beta_2 = 0$, $\beta_3 = s_2^\theta - s_1^\theta = \frac{\Lambda_1^2 - \Lambda_2}{\Lambda_1^3} - \frac{1}{\Lambda_1} = -\frac{\Lambda_2}{\Lambda_1^3}$, and $m_2 = \Lambda_1^2$ implies that

$$0 = 2m_2\beta_1\beta_3 + m_2\left(\beta_2\right)^2 + 3m_3\left(\beta_1\right)^2\beta_2 + m_4\left(\beta_1\right)^4$$
$$= -\frac{2\Lambda_1^2}{\Lambda_1}\frac{\Lambda_2}{\Lambda_1^3} + \frac{m_4}{\Lambda_1^4}.$$

Therefore, we have that $m_4 = 2\Lambda_1^2\Lambda_2$, again verifying (2) for $s = 2$. Proceeding in a similar fashion, all odd-order asymptotic spectral moments can be shown to be zero; while the closed form expressions of all even moments up to the tenth moment are as follows:

$$m_6 = 2\Lambda_1^3\Lambda_3 + 3\Lambda_1^2\Lambda_2^2, \tag{37}$$
$$m_8 = 2\Lambda_1^4\Lambda_4 + 8\Lambda_1^3\Lambda_2\Lambda_3 + 4\Lambda_1^2\Lambda_3,$$
$$m_{10} = 2\Lambda_1^5\Lambda_5 + 10\Lambda_1^4\Lambda_2\Lambda_4 + 5\Lambda_1^4\Lambda_3^2 + 20\Lambda_1^3\Lambda_2^2\Lambda_3 + 5\Lambda_1^2\Lambda_4,$$

which coincide with the closed-form expression provided by our main result (Theorem 1). We would like to remark that the expressions obtained for higher-order values of $s_k^\theta$ and $m_k$ become forbiddingly complex as $k$ increases. In fact, the algebraic manipulations required to compute these values quickly become unmanageable as $k$ increases. In contrast, Theorem 1 provides explicit closed-form expressions for the moments of $F$ (i.e., the asymptotic spectral moments of $\mathbf{A}_n$), which can be used to efficiently compute a truncated sequence of (asymptotic) spectral moments of the random matrix $\mathbf{A}_n$.

## Appendix C. Relations to the band matrix model

The rank one pattern of variances that we specify for the entries of the random matrix ensemble in this work relates to the band matrix model as studied by Anderson and Zeitouni [8]. The band matrix model refers to ensembles of symmetric or Hermitian random matrices with independent entries where variance of each $i, j$-th entry is specified as $\frac{1}{n}f(\frac{i}{n}, \frac{j}{n})$ for a sufficiently well-behaved positive-real symmetric function $f(\cdot,\cdot)\colon [0,1]^2 \to \mathbb{R}^+$. Earlier results on this model consider special cases such as Gaussian entries [33] or with specific requirements on function $f$; for instance Reference [28] shows that with $f(x,y) = \nu(|x-y|)$ for some bounded $\nu(\cdot)$ satisfying $\int_0^1 \nu(t)dt = 1$, the LSD of the band matrix model is the same as Wigner's semi-circle law. This result is recovered under the more general condition $\int_0^1 f(x,y)dy = 1$ in Theorem 3.5 of [8] and in the case of our rank-one variance model we recover exactly the Wigner ensemble. This is because with $f(\frac{i}{n}, \frac{j}{n}) = \sigma_i\sigma_j, \forall i,j$ we can write the integral condition as the limit of a Riemann sum to get that $\int_0^1 f(x,y)dy = \lim_{n\to\infty} \frac{1}{n}\sum_{j=1}^n f(\frac{i}{n}, \frac{j}{n}) = 1$ for some $i$ fixed and now from $f(\frac{i}{n}, \frac{j}{n}) = \sigma_i\sigma_j$ we get that $\sigma_i = \frac{1}{\Lambda_1}, \forall i$, implying $\sigma_i = \sigma_j = 1, \forall i,j$ which retrieves the classical Wigner ensemble described above. In [24] the authors consider the limiting spectral distribution for the sample covariance matrix of a band matrix model; they rely on the Stieltjes tranformation techniques which are an alternative to the moments method that we use in this paper but come short of an explicit characterization as was the case with the free probability and $S$-transform techniques described in Appendix B.

The specification of variances in [8] are in terms of a Polish space of colors denoted by $\mathcal{C}$ which acts as an auxiliary space: any two variables whose indexed locations have the same colors will satisfy the same constraints on their moments. Accordingly, a coloring function $\kappa_0(\cdot)\colon \mathbb{N} \to \mathcal{C}$ maps any index (or letter in the nomenclature of Anderson and Zeitouni [8]), $i \in \mathbb{N}$, to a color $\kappa_0(i) \in \mathcal{C}$. Let $\theta\{\cdot\}$ be the probability distribution on $\mathcal{C}$ that measures the frequency of colors as induced by the mapping $\kappa_0(\cdot)$. Formally $\theta\{\cdot\}$ can be defined as the weak limit as $n \to \infty$ of the measures $\theta_n\{\cdot\} = \frac{1}{n}\sum_{i=1}^n \delta_{\kappa_0(i)}\{\cdot\}$; and indeed, once we specialize the band matrix model of Anderson and Zeitouni [8] to the case of our random matrix ensemble with rank-one variance profile, we shall see $\theta\{\cdot\}$ plays the role of the limiting distribution for the entries of $\{\sigma_i, i \in \mathbb{N}\}$ sequence as was the case in our formulation of the free multiplicative convolution technique in Appendix B. In notation of [8] all moments of entries of the random matrix $\mathbf{A}_n$ are constrained as follows: $\mathbb{E}\left\{|\mathbf{a}_{ii}|^k\right\} \leq d^{(k)}(\kappa_0(i))$, and $\mathbb{E}\left\{|\mathbf{a}_{ij}|^k\right\} \leq s^{(k)}(\kappa_0(i), \kappa_0(j))$ for $i \neq j$, with inequalities being strict when $k = 2, 4$. Here, $d^{(k)}(\cdot)$ and $s^{(k)}(\cdot, \cdot)$ are bounded positive real functions defined on $\mathcal{C}$ and $\mathcal{C} \times \mathcal{C}$ respectively, and they are assumed measurable and almost surely continuous with respect to $\theta\{\cdot\}$. If the preceding moments conditions replace Assumptions 1 to 5 in Section 1, then the main result of [8] concerning the LSD $F(\cdot)$ is stated as follows.

**Proposition 2.** (Theorem 3.2 of [8].) *As $n \to \infty$ the ESD $\mathbf{F}_n(\cdot)$ converge weakly in probability to the distribution $F(\cdot)$ which uniquely satisfies*

$$m_s = \int_{-\infty}^{+\infty} x^s \, dF(x) = \int_{c \in \mathcal{C}} \Phi_{s+1,\theta}(c) d\theta(c), \forall s \in \mathbb{N}_0. \tag{38a}$$

*Here, for all $c \in \mathcal{C}$, $\{\Phi_{s,\theta}(c), s \in \mathbb{N}\}$ is the unique real sequence satisfying*

$$\Phi_\theta(c,t) = \sum_{s=1}^{\infty} \Phi_{s,\theta}(c) t^s, \tag{39a}$$

*and the latter is the formal power series characterized by the generating function identity*

$$\Phi_\theta(c,t) = t \left(1 - t \int_{c' \in \mathcal{C}} s^{(2)}(c,c') \Phi_\theta(c',t) d\theta(c') \right)^{-1}. \tag{40a}$$

When the integral equation (40a) is expanded in terms of the formal series in (39a), it leads to a sequence of recursive equations with the unknowns $\Phi_{s,\theta}(c)$, $s \in \mathbb{N}$ whose solutions upon replacement in (38a) yield the same asymptotic moments as claimed by Theorem 1.

On the one hand, our rank one pattern includes some variance profiles that are not covered by the specification of variances in the color space under the band matrix model. To see how, consider the mapping of indexes to colors under the coloring function $\kappa_0(\cdot)$, and suppose without any loss in generality that the indexes $i \in \mathbb{N}$ are ordered such that indexes with the same color are consecutive to each other. Subsequently, the specification of second moments under $s^2(\cdot,\cdot)$ are such that in any block whose rows have the same color and so do its columns, all entries have the same variance. However, such a configuration excludes many possibilities. In particular, while for separable functions $s^2(\cdot,\cdot)$ any specification of variances under the coloring scheme of the band matrix model can be equivalently captured by a rank one pattern of variances and through a particular sequence $\{\sigma_i, i \in \mathbb{N}\}$; there are certain ways in which the variances may be specified as a rank one pattern, but not by a symmetric separable function $s^2(\cdot,\cdot)$ over the Cartesian product of two copies of a Polish space of colors. As an example, suppose that the sequence $\sigma_i, i \in \mathbb{N}$ is derived by uniform sampling from an exponential function: $\Delta_n e^{-i\delta_n}, i \in [n]$ for some appropriate choice of constants $\Delta_n$ and $\delta_n$. Next note that as $n \to \infty$, those entries $\mathbf{a}_{ij}$ whose sum of row and column indexes are the same, will posses the same asymptotic variances; whereas such specification of variances along the level sets of $i+j$ for $i,j \in \mathbb{N}$ cannot be realized by any assignments of blocks in a Polish color space.

On the other hand, Theorem 1 considers a special case of the band matrix model where the function $s^2(\cdot,\cdot)$, representing the specification of the second order moments on the color space, is separable with the following factorization: $s^2(c,c') = \rho(c)\rho(c')$, $\forall (c,c') \in \mathcal{C} \times \mathcal{C}$ for some appropriately defined function $\rho(\cdot)$. In this particular case which is captured by the rank-one pattern of variances, Theorem 1 provides explicit solutions for the moments satisfying the system of integral equations described in (38a) to (40a) of Proposition 2. In what follows we derive the asymptotic moments of the band matrix model given a rank-one pattern variance and using Proposition 2, as a verification and in agreement with Theorem 1.

First note that per Assumption 3, $\mathbb{E}\left\{\mathbf{a}_{ij}^2\right\} = s^{(2)}(\kappa_0(i),\kappa_0(j)) = \sigma_i\sigma_j$. Let function $\rho\colon \mathcal{C} \to \mathbb{R}^+$ be such that for all $i \in \mathbb{N}$, $\rho(\kappa_0(i)) = \sigma_i$, and in particular we get that $s^{(2)}(c,c') = \rho(c)\rho(c')$ for all $(c,c') \in \mathcal{C}\times\mathcal{C}$. Furthermore, given the probability measures $\theta_n\{\cdot\} = \frac{1}{n}\sum_{i=1}^n \delta_{\kappa_0(i)}\{\cdot\}$ and their weak limit $\theta(\cdot)$, we have that for any $k \in \mathbb{N}$,

$$\int_{c' \in \mathcal{C}} \rho(c')^k d\theta(c') = \lim_{n \to \infty} \int_{c' \in \mathcal{C}} \rho(c')^k d\theta_n(c') \tag{41a}$$

$$= \lim_{n \to \infty} \frac{1}{n} \sum_{i=1}^n \rho(\kappa_0(i))^k$$

$$= \lim_{n \to \infty} \frac{1}{n} \sum_{i=1}^n \sigma_i^k = \Lambda_k.$$

To derive asymptotic relations for the moments per (38a), we begin by replacing $\Phi_\theta(c,t)$ in (40a) with its formal power series expansion given by (39a), and then expanding both sides of (40a) to get:

$$\sum_{s=1}^\infty \Phi_{s,\theta}(c) t^s = t\left(1 - \rho(c)\sum_{s=1}^\infty t^{s+1} \int_{c' \in \mathcal{C}} \rho(c')\Phi_{s,\theta}(c')d\theta(c')\right)^{-1}$$

$$= t + \rho(c)\sum_{s=1}^\infty t^{s+2} \int_{c' \in \mathcal{C}} \rho(c')\Phi_{s,\theta}(c')d\theta(c')$$

$$+ t \sum_{k=2}^\infty \left(\rho(c)\sum_{s=1}^\infty t^{s+1} \int_{c' \in \mathcal{C}} \rho(c')\Phi_{s,\theta}(c')d\theta(c')\right)^k.$$

Equating the coefficients of $t^s$ for each $s \in \mathbb{N}$ yields the following series of equations that are solved for $\Phi_{s,\theta}(c)$, $s \in \mathbb{N}$ recursively. Beginning from $\Phi_{1,\theta}(c) = 1$ and $\Phi_{2,\theta}(c) = 0$, we can proceed in a recursive fashion, while using (41a) to substitute for the respective integrals as they appear. Subsequently, we get

$$\Phi_{3,\theta}(c) = \rho(c)\int_{c'\in\mathcal{C}}\rho(c')\Phi_{1,\theta}(c')d\theta(c') = \rho(c)\int_{c'\in\mathcal{C}}\rho(c')d\theta(c') = \rho(c)\Lambda_1,$$

$$\Phi_{4,\theta}(c) = \rho(c)\int_{c'\in\mathcal{C}}\rho(c')\Phi_{2,\theta}(c')d\theta(c') = 0, \qquad (43a)$$

$$\Phi_{5,\theta}(c) = \rho(c)\int_{c'\in\mathcal{C}}\rho(c')\Phi_{3,\theta}(c')d\theta(c') + \rho(c)^2\left(\int_{c'\in\mathcal{C}}\rho(c')\Phi_{1,\theta}(c')d\theta(c')\right)^2$$

$$= \rho(c)\left(\int_{c'\in\mathcal{C}}\rho(c')^2 d\theta(c')\right)\int_{c'\in\mathcal{C}}\rho(c')d\theta(c') + \rho(c)^2\left(\int_{c'\in\mathcal{C}}\rho(c')d\theta(c')\right)^2$$

$$= \rho(c)\Lambda_2\Lambda_1 + \rho(c)^2\Lambda_1^2,$$

and so forth. By the same token, for general $s \geq 5$ we get

$$\Phi_{s,\theta}(c) = \sum_{\substack{k\in[\lfloor s/2\rfloor] \\ (r_1,\ldots,r_k)\in\mathcal{R}'_{k,s}}} \rho(c)^k \int_{c'\in\mathcal{C}}\rho(c')\Phi_{r_1,\theta}d\theta(c')\ldots\int_{c'\in\mathcal{C}}\rho(c')\Phi_{r_k,\theta}d\theta(c'), \quad (44a)$$

where for all integers $s$ and $1 \leq k \leq \lfloor s/2 \rfloor$, we have defined

$$\mathcal{R}'_{k,s} = \{(r_1,\ldots,r_k) \in \mathbb{N}^k : k + \sum_{j=1}^k r_j = s - 1\}.$$

In particular, note that whenever $s$ is even then for any $1 \leq k \leq \lfloor s/2 \rfloor$ and all $(r_1,\ldots,r_k) \in \mathcal{R}'_{k,s}$ there is some $r_j$, $j \in [k]$ such that $r_j \in [s-2]$ is an even number. To see why consider two cases: if $k$ is odd then $s - 1 - k = \sum_{j=1}^k r_j$ is an even number, and sum of an odd number, $k$, of odd numbers $r_j$, $j \in [k]$ can never be an even number, subsequently at least one of the numbers $r_1$, ..., $r_k$ has to be even or else their sum $\sum_{j=1}^k r_j$ cannot be an even number. Similarly, if $k$ is even, then $s - 1 - k = \sum_{j=1}^k r_j$ is an odd number and the sum of an even number, $k$, of all odd numbers $r_j$, $j \in [k]$ can never be an odd number; hence again, for any $1 \leq k \leq \lfloor s/2 \rfloor$ and all $(r_1,\ldots,r_k) \in \mathcal{R}'_{k,s}$ there is some $r_j$, $j \in [k]$ such that $r_j \in [s-2]$ is an even number.

So far, we have shown that if $s$ is an even number, then for every product term appearing in the summation on the left hand side of (44a), there is some $r_j$, $j \in [k]$ such that $r_j \in [s-2]$ is an even number. This together with the fact that $\Phi_{2,\theta}(c) = \Phi_{4,\theta}(c) \equiv 0$, per above and (43a), imply that $\Phi_{s,\theta}(c) \equiv 0$ for all even integers $s$. Subsequently, we can apply (38a) of Proposition 2 to get that $m_s = \int_{c\in\mathcal{C}}\Phi_{s+1,\theta}(c)d\theta(c)$ for all odd integers $s$, verifying the asymptotic expressions in Theorem 1 for odd moments. For even moments as well we can apply (38a) to the derived expressions in (43a) to get that

$$m_2 = \int_{c \in \mathcal{C}} \Phi_{3,\theta}(c) \Phi_{3,\theta}(c) = \Lambda_1 \int_{c' \in \mathcal{C}} \rho(c') d\theta(c') = \Lambda_1^2,$$
$$m_4 = \int_{c \in \mathcal{C}} \Phi_{5,\theta}(c) \Phi_{3,\theta}(c) = \Lambda_2 \Lambda_1 \int_{c' \in \mathcal{C}} \rho(c') d\theta(c') + \Lambda_1^2 \int_{c' \in \mathcal{C}} \rho(c')^2 d\theta(c') = 2\Lambda_1^2 \Lambda_2.$$

Proceeding in a similar fashion, the closed form expressions of all even moments up to the tenth moment would agree with (37) and correspond exactly to those claimed in (2) by Theorem 1. In conclusion, the asymptotic expressions of the moments in in (2) are giving explicit solutions that would follow if one were to solve the integral equation in (40a), or equivalently, the infinite sequence of recursive equations in (44a).

### Appendix D. Moments method for almost sure and weak convergence

This appendix includes a detailed execution of the analytical steps leading to the conclusions spelled out in Theorem 2 as the method of moments. We begin by establishing the existence and finiteness of the moments in (3). Note that Assumption 2 (boundedness), together with the Gershgorin disk theorem (see Section VIII.6 of [12]), imply that with probability one, $\mathbf{m}_k^{(n)} < n^k K^k$ for all $n, k \in \mathbb{N}$, i.e., all spectral moments are finite. The following lemma is a specialization of Corollary 2.3.6 of [37] to the case of our random matrix ensemble $\mathbf{A}_n(\Psi)$.

**Lemma 7.** (Upper tail estimate for the operator norm.) *There exists absolute, not depending on $n$, constants $C, c > 0$ such that $\mathbb{P}\{\rho(\mathbf{A}_n) > K\gamma\} \leq Ce^{-c\gamma n}$, for all $\gamma \geq C$ and any $n \in \mathbb{N}$.*

**Corollary 3.** (All finite expected spectral moments.) *The expected spectral moments can be bounded as $\bar{m}_{2k}^{(n)} \leq (KC)^{2k} + (2k)!/(cn)^{2k} \leq (KC)^{2k} + (2k)!/c^{2k}$, $\forall k, n \in \mathbb{N}$, where $C, c > 0$ are absolute constants, not depending on $n$. In particular, the expected spectral distributions $\overline{F}_n(\cdot), n \in \mathbb{N}$ have all their moments finite.*

*Proof.* It follows from Lemma 7 that

$$\mathbb{P}\{|\lambda(\mathbf{A}_n)| > \gamma\} \leq \mathbb{P}\{\rho(\mathbf{A}_n) > \gamma\} \leq C\exp(-c\gamma),$$

for all $n \in \mathbb{N}$ and any $\gamma > C$. The claimed bound now follows as

$$\bar{m}_{2k}^{(n)} = \mathbb{E}\left\{|\lambda(\mathbf{A}_n)|^{2k}\right\} = \int_{\alpha=0}^{+\infty} \mathbb{P}\left\{|\lambda(\mathbf{A}_n)|^{2k} > \alpha\right\} d\alpha =$$

$$\int_{\alpha=0}^{KC} \mathbb{P}\left\{|\lambda(\mathbf{A}_n)|^{2k} > \alpha\right\} d\alpha + \int_{\alpha=KC}^{+\infty} \mathbb{P}\left\{|\lambda(\mathbf{A}_n)|^{2k} > \alpha\right\} d\alpha \le$$

$$(KC)^{2k} + \int_{\gamma=0}^{+\infty} Ce^{-c\gamma n} 2k\gamma^{2k-1} d\gamma = (KC)^{2k} + \frac{2k!}{(cn)^{2k}},$$

where we invoked Lemma 7 with the change of variable $\alpha = \gamma^{2k}$ in writing the inequality. □

The preceding results are an instance of the miscellany of results on the concentration of eigenvalues for random matrix ensembles [40]. On this topic, the classical result of Füredi and Komlós [23] gives an almost sure upper bound of $2\sigma\sqrt{n} + O(n^{1/3}\log(n))$, later improved to $O(n^{1/4}\log(n))$ by Vu [48], for the operator norm of a random matrix ensemble with independent entries and identical variances $\sigma^2$. Subsequent results [6, 26] employ the powerful machinery of isoperimetric inequalities for product spaces due to Talagrand [36], and derive sub-exponential bounds for the concentration of the norms and eigenvalues of random matrices around their mean values. In the same venue, Lemma 7, together with the Borel-Cantelli lemma (see Theorem 4.3 of [13]), implies that with probability one, the support of the ESD $\mathbf{F}_n(\cdot)$ is asymptotically compact.

**Corollary 4.** (*Asymptotically compact support.*) *Almost surely, it holds true that* $\limsup_{n\to\infty} \rho(\mathbf{A}_n) < Z$, *for some absolute constant* $Z > 0$, *and in particular we take* $Z > max\{1, 4\Lambda_1\}$.

In particular, we have that almost surely $\rho(\mathbf{A}_n) = O(1)$, which justifies the normalization factor of $\frac{1}{\sqrt{n}}$ used in the definition of the random matrix ensemble $\mathbf{A}_n\left(\{\sigma_i : i \in [n]\}\right)$. Indeed, the local behavior of eigenvalues at the edge of the spectrum is of much interest in applied areas such as quantum theory and statistical mechanics (see Chapter 1 of [27] and [22]). It is the well-known result of Tracy and Widom, who establish the joint distribution of the $k$-largest (or $k$-smallest) eigenvalues of a random matrix with Gaussian entries [42, 43], and the distribution that they derive is later shown, by Soshnikov [34], to apply equally well to the larger class of symmetric Wigner ensembles with independent zero-mean entries and identical variances; a phenomenon known as edge universality [38].

In what follows we prove that the (random) ESD $\mathbf{F}_n(\cdot)$ converges almost surely to the same weak limit as that of the expected spectral distributions $\overline{F}_n(\cdot)$ defined in Subsection 2.1. This is a strong law of large numbers type result that concerns random probability measures endowed with the topology of weak convergence. We thus pave the way from the pointwise convergence of expected spectral moments, $\lim_{n\to\infty} \bar{m}_k^{(n)} = m_k$, $\forall k \in \mathbb{N}$, to the almost sure and weak convergence of the ESD $\mathbf{F}_n(\cdot)$ to the unique distribution $F(\cdot)$ satisfying $\forall k \in \mathbb{N}$, $\int_{-\infty}^{+\infty} x^k \, dF(x) = m_k$.

**Almost sure weak convergence of the ESD**

The first set of results in this subsection allows us to conclude weak convergence of the expected spectral distributions $\overline{F}_n(\cdot)$ to the distribution $F(\cdot)$ as $n \to \infty$ from the pointwise convergence of their moments sequence $\bar{m}_k^{(n)}, k \in \mathbb{N}$ to the sequence of moments $m_k, k \in \mathbb{N}$ defined earlier in Section 1. To this end, we need the following four lemmas, of which the first two are restatements of Theorems 4.5.2 and 4.5.5 in [16], the third is the celebrated test due to Carleman [5] and is proved as Lemma B.3 in [9], and the fourth is a well-known consequence of Helly's selection theorem that establishes the relative compactness of a tight family of probability measures with respect to the topology of weak convergence, which can be found for instance as Theorem 25.10 of [13]. Recall, apropos, that a sequence of distributions $\{\hat{F}_n(\cdot) \colon n \in \mathbb{N}\}$ is said to be tight if for each $\epsilon > 0$ there exist real numbers $x$ and $y$ that $\hat{F}_n(x) < \epsilon$ and $\hat{F}_n(y) > 1 - \epsilon$ for all n.

**Lemma 8.** (Uniform integrability.) *Suppose that $\{\hat{F}_n(\cdot) \colon n \in \mathbb{N}\}$ is a sequence of distribution functions and $\hat{F}(\cdot)$ is a distribution function, such that $\hat{F}_n(\cdot)$ converges weakly to $\hat{F}(\cdot)$ as $n \to \infty$. Suppose further that for some $s \in \mathbb{N}$ and $M > 0$, the even finite moments given by $\hat{m}_{2s}^{(n)} = \int_{-\infty}^{+\infty} x^{2s} d\hat{F}_n(x)$ satisfy $\sup_{n \in \mathbb{N}} \hat{m}_{2s}^{(n)} < M$. Then $\forall k \in [2s - 1]$, $\lim_{n\to\infty} \hat{m}_k^{(n)} = \int_{-\infty}^{+\infty} x^k d\hat{F}(x)$.*

**Lemma 9.** (The method of moments.) *Suppose that there is a unique distribution function $\hat{F}(\cdot)$ associated with the sequence of moments $\{\hat{m}_k \colon k \in \mathbb{N}\}$, all finite; such that $\hat{m}_k = \int_{-\infty}^{+\infty} x^k d\hat{F}(x)$, $\forall k \in \mathbb{N}$. Suppose further that for all $n \in \mathbb{N}$, $\hat{F}_n(\cdot)$ is a distribution function, which has all its moments finite and given by $\hat{m}_k^{(n)} = \int_{-\infty}^{+\infty} x^k d\hat{F}_n(x)$, $\forall k \in \mathbb{N}$. Finally, suppose that for every $k \geq 1$, $\hat{m}_k^{(n)} \to \hat{m}_k$ as $n \to \infty$. Then $\hat{F}_n(\cdot)$* *converges weakly to $\hat{F}(\cdot)$ as $n \to \infty$.*

**Lemma 10.** (Carleman's criterion.) *Suppose that there is a distribution function $\hat{F}(\cdot)$ associated with the sequence of moments $\{\hat{m}_k \colon k \in \mathbb{N}\}$, all finite, such that $\hat{m}_k = \int_{-\infty}^{+\infty} x^k d\hat{F}(x), \forall k \in \mathbb{N}$. Suppose further that $\sum_{k=0}^{\infty} \hat{m}_{2k}^{-1/2k} = \infty$. Then $\hat{F}(\cdot)$ is the unique distribution satisfying $\hat{m}_k = \int_{-\infty}^{+\infty} x^k d\hat{F}(x)$, for all $k \in \mathbb{N}$.*

**Lemma 11.** (Helly's selection principle.) *Given a sequence of distributions $\{\hat{F}_n(\cdot) \colon n \in \mathbb{N}\}$, its tightness is a necessary and sufficient condition that for every subsequence $\hat{F}_{n_k}(\cdot)$, $\{n_k \colon k \in \mathbb{N}\} \subset \mathbb{N}$, there exist a further subsequence $\hat{F}_{n_{k_j}}(\cdot)$, $\{n_{k_j} \colon j \in \mathbb{N}\} \subset \{n_k \colon k \in \mathbb{N}\}$, and a distribution $\hat{F}(\cdot)$, such that $\hat{F}_{n_{k_j}}(\cdot)$ converges weakly to $\hat{F}(\cdot)$ as $j \to \infty$.*

We now have all the necessary tools at our disposal to conclude the weak convergence of expected spectral distributions from the pointwise convergence of their moments, and as well as to conclude the almost sure weak convergence of the ESD from the almost sure pointwise convergence of the spectral moments; the facts of which are established by the proceeding theorem and the subsequent corollary.

**Theorem 6.** (Existence and uniqueness of LSD.) *If $\forall k \in \mathbb{N}$, $\lim_{n \to \infty} \bar{m}_k^{(n)} = m_k$, then $\overline{F}_n(\cdot)$ converges weakly to $F(\cdot)$ as $n \to \infty$, where $F(\cdot)$ is the unique distribution function satisfying $\forall k \in \mathbb{N}$, $\int_{-\infty}^{+\infty} x^k dF(x) = m_k$.*

*Proof.* The sequence of distribution functions $\{\overline{F}_n(\cdot) \colon n \in \mathbb{N}\}$ have all their moments finite per Corollary 3. To invoke Lemma 9 for the method of moments then, the gist of the proof is in establishing the tightness property for the sequence $\{\overline{F}_n(\cdot) \colon n \in \mathbb{N}\}$ and then verifying Carleman's criterion for the moments sequence $\{m_k \colon k \in \mathbb{N}\}$. Indeed, from Corollary 3 we get that $\sup_{n \in \mathbb{N}} \bar{m}_2^{(n)} < (KC)^4 + 24/c^4 = M_2$ for some absolute constant $M_2 > 0$, and by the inequality of Chebyshev we get for all $y > 0$ and any $n \in \mathbb{N}$ that $\overline{\mathcal{L}}_n\{|\lambda(\mathbf{A}_n)| > y\} = \overline{F}_n(-y) + 1 - \overline{F}_n(y) < M_2/y^2$. Next, since $\max\{\overline{F}_n(-y), 1 - \overline{F}_n(y)\} \leq \overline{F}_n(-y) + 1 - \overline{F}_n(y) < M_2/y^2$, for any $\epsilon > 0$ and all $n \in \mathbb{N}$, we can set $y = (M_2/\epsilon)^{1/2}$ and $x = -y$ to get that $\overline{F}_n(x) < \epsilon$ and $\overline{F}_n(y) > 1 - \epsilon$, whence follows the tightness of $\{\overline{F}_n \colon n \in \mathbb{N}\}$. This tightness per Lemma 11 implies that for some subsequence $\overline{F}_{n_j}(\cdot)$, $\{n_j \colon j \in \mathbb{N}\} \subset \mathbb{N}$ and a subsequential limit, call it $F(\cdot)$, we have that $\overline{F}_{n_j}(\cdot)$ converges weakly to $F(\cdot)$ as $j \to \infty$. But then it has to be that $\forall k \in \mathbb{N}$, $\int_{-\infty}^{+\infty} x^k \, dF(x) = m_k$. Indeed, take any $k \in \mathbb{N}$ fixed, the finite limit $\lim_{j \to \infty} \bar{m}_{2k}^{(n_j)} = m_{2k}$ implies that for some $M_{2k} > 0$ large enough, $\sup_{n_j:\, j \in \mathbb{N}} \bar{m}_{2k}^{(n_j)} < M_{2k}$; thence by Lemma 8 for $s < 2k$ we get that $m_s = \lim_{j \to \infty} \hat{m}_s^{(n_j)} = \int_{-\infty}^{+\infty} x^s d\, F(x)$. In particular, we have shown the existence of a distribution $F(\cdot)$, satisfying $\int_{-\infty}^{+\infty} x^k d\, F(x) = m_k$ for any positive integer $k$. Now, that $F(\cdot)$ satisfies the preceding uniquely follows per Lemma 10 as the sequence $\{m_k \colon k \in \mathbb{N}\}$ given by (2) passes Carleman's test. Effectually, using the bound $\Lambda_k \leq \Lambda_1^k$, true for all $k \in \mathbb{N}$, we get that for all $s \in \mathbb{N}$

$$m_{2s} = \sum_{(r_1,\ldots,r_s) \in \mathcal{R}_s} \frac{2}{s+1} \binom{s+1}{r_1,\ldots,r_s} \Lambda_1^{r_1} \Lambda_2^{r_2} \ldots \Lambda_s^{r_s} \leq \sum_{(r_1,\ldots,r_s) \in \mathcal{R}_s} 2 \binom{s+1}{r_1,\ldots,r_s} \Lambda_1^{s+1}$$
$$= C_s \Lambda_1^{s+1} \leq (4\Lambda_1)^{s+1}, \tag{47a}$$

where in the penultimate equality we have invoked Corollary 4 and the last inequality is by the fact that $C_s < 4^s$, $\forall s \in \mathbb{N}$ (see, for instance, the proof of Lemma 2.1.7 in [7]). Verification of Carleman's criterion is now immediate as

$$\sum_{k=0}^{\infty} m_{2k}^{-1/2k} \geq \sum_{k=0}^{\infty} (4\Lambda_1)^{-(k+1)/2k} \geq \min\left\{ \sum_{k=0}^{\infty} \frac{1}{2\sqrt{\Lambda_1}}, \sum_{k=0}^{\infty} \frac{1}{4\Lambda_1} \right\} = \infty.$$

We have thus established the existence of a distribution $F(\cdot)$, uniquely satisfying $\int_{-\infty}^{+\infty} x^k \, dF(x) = m_k$, $\forall k \in \mathbb{N}$ so that per Lemma 9 $\lim_{n \to \infty} \bar{m}_k^{(n)} = m_k$, $\forall k \in \mathbb{N}$ implies the weak converge of the corresponding distribution functions $\overline{F}_n(\cdot)$ to the distribution $F(\cdot)$; thence, completing the proof.  □

**Corollary 5.** (Almost sure weak convergence of the ESD.) *If $\lim_{n \to \infty} \mathbf{m}_k^{(n)} = m_k$, almost surely for all $k \in \mathbb{N}$, then with probability one $\mathbf{F}_n(\cdot)$ converges weakly to $F(\cdot)$ as $n \to \infty$, where $F(\cdot)$ is the unique distribution function satisfying $\forall k \in \mathbb{N}$, $\int_{-\infty}^{+\infty} x^k \, dF(x) = m_k$.*

*Proof.* By the countable intersection of full probability measure sets, we get that restricted to a measurable set $\mathcal{S} \subset \Omega$ with $\mathbb{P}\{\mathcal{S}\} = 1$, the sequence $\mathbf{F}_n(\cdot)$ have all finite moments, satisfying $\mathbf{m}_2^{(n)} < n^2 K^2$ and $\lim_{n \to \infty} \mathbf{m}_2^{(n)} = m_2$ so that restricted to $\mathcal{S}$ and for some $\mathbf{M}$ not dependent on $n$, we have $\mathbf{m}_2^{(n)} < \mathbf{M}$. Thence, the proof of Theorem 6 applies equally well to the sequence $\mathbf{F}_n(\cdot)$ on this full probability measure set $\mathcal{S}$ and the conclusion of the corollary follows, per Definition 1 of Appendix A.  □

**Weak almost sure convergence of the ESD**

The key in establishing the almost sure convergence of the ESD $\mathbf{F}_n(\cdot)$ to the deterministic weak limit $F(\cdot)$ is in verifying that $\lim_{n\to\infty} \mathbf{m}_k^{(n)} = m_k$, almost surely, for each $k \in \mathbb{N}$. This is achieved through an application of Talagrand's concentration inequality restated from Theorem 2.1.13 of [37] as Lemma 12 below, followed by the Borel-Cantelli lemma leading to the claimed almost sure convergence.

**Lemma 12.** (Talagrand's concentration inequality.) *For each $n \in \mathbb{N}$, let $f_n(\cdot) : \mathbb{R}^{n(n+1)/2} \to \mathbb{R}$ be a convex function acting on the diagonal and upper diagonal entries of the random matrix $\mathbf{A}_n$. Further let $f_n(\cdot)$ be 1-Lipschitz with respect to the Euclidean norm on $\mathbb{R}^{n(n+1)/2}$. Then for any $\lambda$ one has that $\mathbb{P}\{|f(\mathbf{A}_n) - \mathbb{E}\{f(\mathbf{A}_n)\}| \geq \lambda K\} \leq C e^{-c\lambda^2}$, for some absolute, not depending on $n$, constants $C, c > 0$.*

**Lemma 13.** ($k$-Schatten norms.) *For all $n, k \in \mathbb{N}$, $\|\mathbf{A}_n\|_k = \left(\sum_{i=1}^n |\lambda_i(\mathbf{A}_n^k)|^k\right)^{\frac{1}{k}}$ is a convex function, mapping the diagonal and upper diagonal entries of the random matrix $\mathbf{A}_n$ to positive reals. Furthermore, for $k \geq 2$ it is $\sqrt{2}$-Lipschitz with respect to the Euclidean distance on $\mathbb{R}^{n(n+1)/2}$ and for $k = 1$ it is $\sqrt{2n}$-Lipschitz with respect to that same metric on $\mathbb{R}^{n(n+1)/2}$.*

*Proof.* First note that $\|\mathbf{A}_n\|_k$ is the $k$-Schatten matrix norm for $\mathbf{A}_n$ and is therefore convex (see Section IV.2 of [12]). The Lipschitz property is a consequence of the reverse triangle inequality of the norms. Indeed, for all $k$ and any two $n \times n$ matrices $A$ and $B$, it holds true that $|\|A\|_k - \|B\|_k| \leq \|A - B\|_k$. The $\sqrt{2}$-Lipschitz for $k \geq 2$ property now follows as

$$\|\mathbf{A}_n\|_k^k = \sum_{i=1}^n |\lambda_i(\mathbf{A}_n)|^k \leq \left(\sum_{i=1}^n |\lambda_i(\mathbf{A}_n)|^2\right)^{k/2} = \|\mathbf{A}_n\|_2^k \leq \sqrt{2}\left(\sum_{1\leq i \leq j \leq n} \mathbf{a}_{ij}^2\right)^{k/2},$$

where in the last inequality we used the fact that the 2-Schatten norm $\|\mathbf{A}_n\|_2$ coincides with the Frobenius norm of $\mathbf{A}_n$. Similarly for $k = 1$, the Cauchy-Schwartz inequality implies that

$$\|\mathbf{A}_n\|_1 = \sum_{i=1}^n |\lambda_i(\mathbf{A}_n)| \leq \sqrt{n}\left(\sum_{i=1}^n |\lambda_i(\mathbf{A}_n)|^2\right)^{1/2} \leq \sqrt{2n}\left(\sum_{1\leq i \leq j \leq n} \mathbf{a}_{ij}^2\right)^{1/2},$$

which, together with the reverse triangle inequality, implies the claimed $\sqrt{2n}$-Lipschitz property for $\|\mathbf{A}_n\|_1$. □

**Lemma 14.** (Almost sure convergence of the spectral moments.) *For any $k \in \mathbb{N}$, if $\lim_{n \to \infty} \bar{m}_k^{(n)} = m_k$, then $\lim_{n \to \infty} \mathbf{m}_k^{(n)} = m_k$, almost surely.*

*Proof.* We first prove the claim directly for all even $k$ and then by induction for odd $k$ as well. For $k$ any fixed even integer per Lemma 13, $\frac{1}{\sqrt{2}}\|\mathbf{A}_n\|_k = \frac{1}{\sqrt{2}}\text{trace}(\mathbf{A}_n^k)^{1/k}$ is a convex 1-Lipschitz function satisfying the presumptions of Lemma 12, which together with (3), implies that for any $k, n \in \mathbb{N}$,

$$\bar{C}e^{-\bar{c}\lambda^2} \geq \mathbb{P}\left\{\left|\frac{1}{\sqrt{2}}\text{trace}(\mathbf{A}_n^k)^{1/k} - \mathbb{E}\left\{\frac{1}{\sqrt{2}}\text{trace}(\mathbf{A}_n^k)^{1/k}\right\}\right| \geq \lambda K\right\}$$
$$= \mathbb{P}\left\{\left|(\mathbf{m}_k^{(n)})^{1/k} - (\bar{m}_k^{(n)})^{1/k}\right| \geq \lambda\sqrt{2}K/n^{1/k}\right\},$$

for some absolute, not depending on $n$, constants $\bar{C}, \bar{c} > 0$. Using the change of variables $\epsilon = \lambda\sqrt{2}K/n^{1/k}$ and $\hat{c} = \bar{c}/(2K^2)$, we then get that

$$\bar{C}e^{-\hat{c}\epsilon^2 n^{2/k}} \geq \mathbb{P}\left\{\left|(\mathbf{m}_k^{(n)})^{1/k} - (\bar{m}_k^{(n)})^{1/k}\right| \geq \epsilon\right\}.$$

In particular, we have for any $\epsilon > 0$ that

$$\sum_{n=1}^{\infty} \mathbb{P}\left\{\left|(\mathbf{m}_k^{(n)})^{1/k} - (\bar{m}_k^{(n)})^{1/k}\right| \geq \epsilon/2\right\} < \infty,$$

thence the Borel-Cantelli lemma (Theorem 4.3 of [13]) implies that with probability one $\lim_{n \to \infty} (\mathbf{m}_k^{(n)})^{1/k} = \lim_{n \to \infty} (\bar{m}_k^{(n)})^{1/k}$, and the claim for even $k$ now follows by taking the $k$-th power and invoking continuity. Next consider the case of odd $k$, and note that per Corollary 4 with probability one for $n > \hat{N}$ large enough, $\mathbf{A}_n + ZI_n$ has all its eigenvalues positive so that $\|\mathbf{A}_n + ZI_n\|_k$ is convex and $\sqrt{2}$-Lipschitz for $k > 1$ and $\sqrt{2n}$-Lipschitz for $k = 1$. Indeed, for $k = 1$ Lemma 12 implies that

$$\bar{C}_1 e^{-\bar{c}_1\lambda^2} \geq \mathbb{P}\left\{\left|\frac{1}{\sqrt{2n}}\text{trace}(\mathbf{A}_n + ZI_n) - \mathbb{E}\left\{\frac{1}{\sqrt{2n}}\text{trace}(\mathbf{A}_n + ZI_n)\right\}\right| \geq \lambda K\right\}$$
$$= \mathbb{P}\left\{\frac{1}{n}\left|\text{trace}(\mathbf{A}_n + ZI_n) - \frac{1}{n}\mathbb{E}\left\{\text{trace}(\mathbf{A}_n + ZI_n)\right\}\right| \geq \lambda\sqrt{2}K/\sqrt{n}\right\},$$

and the change of variable $\hat{\epsilon} = \lambda\sqrt{2}K/\sqrt{n}$ again yields

$$\mathbb{P}\left\{\left|\frac{1}{n}\text{trace}(\mathbf{A}_n + ZI_n) - \frac{1}{n}\mathbb{E}\left\{\text{trace}(\mathbf{A}_n + ZI_n)\right\}\right| \geq \hat{\epsilon}\right\} \leq \bar{C}_1 e^{-\bar{c}_1\hat{\epsilon}^2 n/(2K^2)},$$

which is summable over $n$, implying per Borel-Cantelli lemma that $\frac{1}{n}\text{trace}(\mathbf{A}_n + ZI_n)$ converges almost surely to $\frac{1}{n}\mathbb{E}\{\text{trace}(\mathbf{A}_n + ZI_n)\}$ as $n \to \infty$, or equivalently that $\lim_{n\to\infty} \frac{1}{n} \mathbf{m}_1^{(n)} = m_1$ almost surely, completing the proof for $k = 1$. The proof for odd $k > 1$ proceeds by induction for which the base $k = 1$ is already established. Fix $k$ and suppose that for all odd integers less than $k$, $\lim_{n\to\infty} \mathbf{m}_k^{(n)} = m_k$ almost surely. Repeating the above argument for $k$ yields that $\left(\frac{1}{n}\text{trace}(\mathbf{A}_n + ZI_n)^k\right)^{1/k}$ as $n \to \infty$, converges almost surely to

$$\mathbb{E}\left\{\left(\frac{1}{n}\text{trace}(\mathbf{A}_n + ZI_n)^k\right)^{1/k}\right\},$$

implying by continuity that

$$\lim_{n\to\infty} \frac{1}{n} \sum_{i=1}^n (\lambda_i(\mathbf{A}_n) + Z)^k = \lim_{n\to\infty} \frac{1}{n} \mathbb{E}\left\{\sum_{i=1}^n (\lambda_i(\mathbf{A}_n) + Z)^k\right\},$$

almost surely. Expanding the binomial terms and invoking the induction hypothesis, together with the claim proved for the even case, then yields that $\lim_{n\to\infty} \frac{1}{n} \sum_{i=1}^n \lambda_i^k(\mathbf{A}_n) = \lim_{n\to\infty} \frac{1}{n} \mathbb{E}\{\sum_{i=1}^n \lambda_i^k(\mathbf{A}_n)\}$ almost surely, completing the proof by induction for all odd $k$. □

While Corollary 5 gives already the almost sure weak convergence of the ESD as implied by the almost sure convergence of their moments sequence, the path from the latter to the weak almost sure convergence of $\mathbf{F}_n(\cdot)$ is paved by the next proposition, which is the last result in this section, providing us with all that is needed for concluding Theorem 1 from the pointwise convergence of the expected spectral moment sequence.

**Proposition 3.** (Weak almost sure convergence of the ESD.) *If $\lim_{n\to\infty} \mathbf{m}_k^{(n)} = m_k$, with probability one for all $k \in \mathbb{N}$, then $\mathbf{F}_n(\cdot)$ converges weakly, almost surely, to $F(\cdot)$ as $n \to \infty$, where $F(\cdot)$ is the unique distribution function satisfying $\forall k \in \mathbb{N}$, $\int_{-\infty}^{+\infty} x^k \, dF(x) = m_k$.*

*Proof.* For any fixed bounded continuous real function $f(\cdot): \mathbb{R} \to \mathbb{R}$ and some real constant $B > 0$ such that $f(x) < B, \forall x \in \mathbb{R}$, it is required to show that

$$\lim_{n\to\infty} \int_{-\infty}^{+\infty} f(x) \, d\mathbf{F}_n(x) = \int_{-\infty}^{+\infty} f(x) \, dF(x),$$

almost surely, cf. Definition 2 of Appendix A. Take $Z$ as in Corollary 4; by the Stone – Weierstrass theorem (Theorem 7.26 of [32]), for some fixed polynomial $q_\epsilon(x) = \sum_{i=0}^L b_i x^i$ and $\epsilon > 0$ we have $\sup_{|x|\leq Z} |q_\epsilon(x) - f(x)| < \epsilon/2$. Adding and subtracting

40                                                                 Preciado and Rahimianthe terms $\int_{-\infty}^{+\infty} q_\epsilon(x)\, d\mathbf{F}_n(x)$ and $\int_{-\infty}^{+\infty} q_\epsilon(x)\, dF(x)$, and then applying the triangle and Jensen inequalities yields

$$\left|\int_{-\infty}^{+\infty} f(x)\, d\mathbf{F}_n(x) - \int_{-\infty}^{+\infty} f(x)\, dF(x)\right| \leq \epsilon + \int_{\{|x|>Z\}} (B + |q_\epsilon(x)|)\, d\mathbf{F}_n(x)$$
$$+ \left|\int_{-\infty}^{+\infty} q_\epsilon(x)\, d\mathbf{F}_n(x) - \int_{-\infty}^{+\infty} q_\epsilon(x)\, dF(x)\right|$$
$$+ \left|\int_{\{|x|>Z\}} (B + |q_\epsilon(x)|)\, dF(x)\right|, \tag{59a}$$

where we also used the facts that

$$\left|\int_{-\infty}^{+\infty} (f(x) - q_\epsilon(x)) d\mathbf{F}_n(x)\right| \leq \int_{-Z}^{+Z} |f(x) - q_\epsilon(x)|\, d\mathbf{F}_n(x)$$
$$+ \int_{\{|x|>Z\}} (B + |q_\epsilon(x)|)\, d\mathbf{F}_n(x),$$

and similarly that

$$\left|\int_{-\infty}^{+\infty} (f(x) - q_\epsilon(x))\, dF(x)\right| \leq \left|\int_{-Z}^{+Z} (f(x) - q_\epsilon(x))\, dF(x)\right|$$
$$+ \int_{\{|x|>Z\}} (B + |q_\epsilon(x)|)\, dF(x),$$

together with

$$\left|\int_{-Z}^{+Z} (f(x) - q_\epsilon(x))\, d\mathbf{F}_n(x)\right| < \epsilon/2$$

and

$$\left|\int_{-Z}^{+Z} (f(x) - q_\epsilon(x))\, dF(x)\right| < \epsilon/2.$$

The almost sure convergence of the moments implies that

$$\lim_{n\to\infty} \left|\int_{-\infty}^{+\infty} q_\epsilon(x)\, d\mathbf{F}_n(x) - \int_{-\infty}^{+\infty} q_\epsilon(x)\, dF(x)\right| = 0,$$

almost surely. Moreover,

$$\int_{\{|x|>Z\}} |x|^k\, dF(x) \leq \frac{m_{2k}}{Z^k} \leq \frac{(4\Lambda_1)^{k+1}}{Z^k}, \forall k \in \mathbb{N}_0, \tag{65a}$$

having used $1 \leq |x|^k/Z^k$ when $|x| > Z$ for the penultimate inequality and (47a) for the last inequality. Notice that by the choice of $Z > 4\Lambda_1$ in Corollary 4 the right hand

side of (65a) is monotonically decreasing in $k$ and goes to zero as $k \to \infty$, whereas by the choice of $Z > 1$ the left hand side of (65a) is increasing in $k$ and it should therefore be true that $\sup_k \int_{-\infty}^{+\infty} |x|^k \mathbb{1}_{\{|x|>Z\}}(x) \, dF(x) \leq \lim_{k \to \infty} (4\Lambda_1)^{k+1}/Z^k = 0$, thence

$$\int_{-\infty}^{+\infty} |x|^k \mathbb{1}_{\{|x|>Z\}}(x) \, dF(x) = 0, \forall k \in \mathbb{N}_0$$

and in particular $\int_{-\infty}^{+\infty} (B + |q_\epsilon(x)|) \mathbb{1}_{\{|x|>Z\}}(x) \, dF(x) = 0$. By the same token and after invoking the almost sure limit of the spectral moments we get that $\lim_{n \to \infty} \int_{\{|x|>Z\}} (B + |q_\epsilon(x)|) \, d\mathbf{F}_n(x)| = 0$, almost surely. The desired conclusion now follows, after taking the limit as $n \to \infty$ of both sides in (59a) to get that for any $\epsilon > 0$ fixed above

$$\lim_{n \to \infty} \left| \int_{-\infty}^{+\infty} f(x) \, d\mathbf{F}_n(x) - \int_{-\infty}^{+\infty} f(x) \, dF(x) \right| \leq \epsilon,$$

almost surely, and choosing $\epsilon > 0$ arbitrarily small. $\square$

We have thus pieced together all the ingredients required of the conclusion of Theorem 2. If $\lim_{n \to \infty} \bar{m}_k^{(n)} = m_k, \forall k \in \mathbb{N}$, then Lemma 14 implies that $\lim_{n \to \infty} \mathbf{m}_k^{(n)} = m_k$, almost surely for any $k \in \mathbb{N}$. Corollary 5 and Proposition 3 would then imply the weak and almost sure convergence of the ESD $\mathbf{F}_n(\cdot)$ as claimed by Theorem 2.

## Acknowledgements

This work was supported by the National Science Foundation grants CNS-1302222 "NeTS: Medium: Collaborative Research: Optimal Communication for Faster Sensor Network Coordination", and IIS-1447470 "BIGDATA: Spectral Analysis and Control of Evolving Large Scale Networks".